\documentclass[twoside,leqno]{article}

\usepackage[letterpaper]{geometry}

\usepackage{siamproceedings}

\usepackage[T1]{fontenc}
\usepackage{amsfonts}
\usepackage{graphicx}
\usepackage{epstopdf}
\usepackage{enumitem}
\usepackage{algorithmic}
\usepackage{amssymb}

\ifpdf
  \DeclareGraphicsExtensions{.eps,.pdf,.png,.jpg}
\else
  \DeclareGraphicsExtensions{.eps}
\fi

\newsiamremark{remark}{Remark}
\newsiamremark{hypothesis}{Hypothesis}
\crefname{hypothesis}{Hypothesis}{Hypotheses}
\newsiamthm{claim}{Claim}

\usepackage{amsopn}

\usepackage{pdfsync}
\usepackage{esint} 

\newcommand{\C}{{\mathbb C}}       
\newcommand{\R}{{\mathbb R}}       

\newcommand{\bS}{{\mathbb S}}

\newcommand{\DD}{{\mathcal D}}

\newcommand{\HH}{{\mathcal H}}

\newcommand{\cN}{{\mathcal N}}
\newcommand{\RR}{{\mathcal R}}

\newcommand{\EE}{{\mathcal E}}

\newcommand{\diam}{{\rm diam}}
\newcommand{\dist}{{\rm dist}}

\newcommand{\capp}{{\operatorname{Cap}}}
\newcommand{\rf}[1]{{(\ref{#1})}}

\newcommand{\supp}{\operatorname{supp}}

\newcommand{\vphi}{{\varphi}}
\newcommand{\ve}{{\varepsilon}}
\newcommand{\vv}{{\vspace{2mm}}}

\newcommand{\wt}[1]{{\widetilde{#1}}}

\newcommand{\pom}{{\partial\Omega}}

\newcommand{\pv}{\operatorname{pv}}

\newcommand{\rad}{{\operatorname{rad}}}

\def\Xint#1{\mathchoice
{\XXint\displaystyle\textstyle{#1}}%
{\XXint\textstyle\scriptstyle{#1}}%
{\XXint\scriptstyle\scriptscriptstyle{#1}}%
{\XXint\scriptscriptstyle\scriptscriptstyle{#1}}%
\!\int}
\def\XXint#1#2#3{{\setbox0=\hbox{$#1{#2#3}{\int}$ }
\vcenter{\hbox{$#2#3$ }}\kern-.58\wd0}}

\def\avint{\;\Xint-}

\def\BMO{\mathop\mathrm{BMO}} 					
\def\Lip{\mathop\mathrm{Lip}} 	                       
\def\H11{\textup{H}^1_1} 					

\theoremstyle{remark}
\newtheorem{question}[theorem]{Question}

\begin{document}

\title{\Large Interactions between quantitative rectifiability, singular integrals, and boundary value problems for harmonic functions}

\author{Xavier Tolsa\thanks{Instituci\'o Catalana de Recerca i Estudis Avan\c cats (ICREA), Universitat Aut\`onoma de Barcelona, and Centre de Recerca Matem\`atica (\email{xavier.tolsa@uab.cat}).}}

\date{}

\maketitle





\fancyfoot[R]{\scriptsize{Copyright \textcopyright\ 2026\\
Copyright for this paper is retained by the author}}


\pagenumbering{arabic}
\setcounter{page}{1}

\begin{abstract} 
This paper surveys different topics where the theory of quantitative rectifiability plays a central role. First, it reviews the characterization
of rectifiability in terms of square functions involving $\beta$ type coefficients and the $\ve^2$ conjecture of Carleson. 
It also discusses the deep connections between rectifiability and the $L^2$ boundedness of Riesz transforms and their application to the Painlev\'e problem for Lipschitz harmonic functions.
Finally, the paper explores recent major advances in connection with harmonic measure and the $L^p$ solvability of the Dirichlet, regularity, and Neumann
problems for the Laplace equation in rough domains, emphasizing the key role of quantitative rectifiability in these developments.
\end{abstract}

\section{Introduction.}

The study of rectifiable sets and measures is one of the main objectives of geometric measure theory.
Recall that a set $E\subset \R^d$ is called $n$-rectifiable, for some natural number $n\leq d$, if there is a countable collection of Lipschitz maps $g_i:\R^n\to\R^d$
such that
$$\HH^n\Big(E\setminus \bigcup_i g_i(\R^n)\Big)=0,$$
where $\HH^n$ stands for the $n$-dimensional Hausdorff measure.
On the other hand, a Radon measure $\mu$ is $n$-rectifiable if it is of the form $\mu = g\,\HH^n|_E$, where $E\subset\R^d$ is 
$n$-rectifiable and $g\in L^1_{loc}(\HH^n|_E)$.
The notion of rectifiability  plays an important role in many problems in analysis, partial differential equations, and geometry.
In the 1990's, very influential works of Jones \cite{JonesTraveling} and David and Semmes \cite{DavidSemmes2} introduced techniques from harmonic analysis to study the $L^2$ boundedness of singular integral operators over rectifiable sets giving rise to what now is known as quantitative rectifiability.

The theory of quantitative rectifiability has been very successful. Perhaps the first spectacular application was in the solution by David \cite{David-vitus} of Vitushkin's conjecture about removable singularities for bounded analytic functions in 1998, for sets with finite length. In the last few years, the techniques of quantitative rectifiability have
found many other outstanding applications to questions in connection with harmonic measure, boundary value problems for elliptic PDE's, or free boundary problems. In this paper we survey some of these advances. There is no attempt of completeness and I will focus mainly on the works that I know better and that are connected to my own research. In fact, this area has broadened so much that writing an exposition covering all the relevant related results would take us too far.

The main topics we will discuss are the following. In Section \ref{sec3} we will recall the traveling salesman theorem of Peter Jones and the notion of uniform rectifiability of David and Semmes. We will also review the characterization of uniformly rectifiable sets 
and measures in terms of the so called $\beta_p$ coefficients and other related square functions and we will see the extension of these 
results to other more general sets in $\R^n$, in particular by Azzam and the author \cite{AzzamTolsa}, \cite{Tolsa-betas}, and Edelen, Naber, and Valtorta \cite{EdelenNaberValtorta}, \cite{NaberValtorta-Rectifiable}. Further, we will discuss the more recent solution of Carleson's $\ve^2$ conjecture about tangent point of
Jordan curves by Jaye, the author, and Villa \cite{JayeTolsaVilla}, and the extension to higher dimensions by Fleschler, the author, and Villa \cite{FleschlerTolsaVilla-Inventiones}.

In Section \ref{sec4} we will describe the connection between the $L^2$ boundedness of singular integral operators, such as Riesz transforms, with rectifiability and uniform rectifiability. Among other topics, we will review the solution of the David-Semmes problem by Nazarov, the author, and Volberg \cite{NazarovTolsaVolberg1}, and other more recent results in collaboration with D\k{a}browski \cite{DabrowskiTolsa}. In Section \ref{sec5} we will explain the applications to the Painlev\'e problem for bounded analytic functions and for Lipschitz harmonic functions in higher dimensions.

Section \ref{sec6} deals with some applications of the theory of quantitative rectifiability to questions concerning harmonic measure.
In particular, we will discuss 
 the one-phase and two-phase problems for harmonic measure, where the Riesz transform plays a fundamental role.  
Section \ref{sec7} is devoted to 
boundary value problems for harmonic functions in rough domains. 
 We will  focus on the $L^p$ solvability of Dirichlet, regularity, and
Neumann problems for the Laplace equation in  domains with uniformly rectifiable boundaries. This is a topic essentially initiated with
the fundamental work of Dahlberg \cite{Dahlberg}, who proved the $L^2$ solvability of the Dirichlet problem for the Laplacian in Lipschitz domains in the late 1970's. 


\section{Preliminaries.}\label{secprelim}

In the paper, constants denoted by $C$ or $c$ depend just on the dimension unless otherwise stated. We will write $a\lesssim b$ if there is $C>0$ such that $a\leq Cb$. We write $a\approx b$ if $a\lesssim b\lesssim a$.
 
 A domain is a connected open set.
 Open balls in $\R^{d}$ centered at $x$ with radius $r>0$ are denoted by $B(x,r)$, and closed balls by 
$\bar B(x,r)$. For an open or closed ball $B\subset\R^{d}$ with radius $r$, we write $\rad(B)=r$.
We use the two notations $S(x,r)\equiv \partial B(x,r)$ for spheres
in $\R^{d}$ centered at $x$ with radius $r$, so that $\bS^{d-1}=S(0,1)$.

All measures are assumed to be Borel measures. The Hausdorff $n$-dimensional measure $\HH^n$ is defined as follows. For $E\subset \R^d$ and $\ve\in (0,\infty]$,
we set
$$\HH^n_\ve(E) = c_n \inf\big\{\sum_i \diam(A_i)^n:\,E\subset\bigcup_i A_i,\,\diam(A_i)\leq \ve\big\},$$
and 
$$\HH^n(E) = \sup_{\ve>0}\HH^n_\ve(E) = \lim_{\ve\to0}\HH^n_\ve(E),$$
where the constant $c_n$ is chosen so that $\HH^n$ coincides with the Lebesgue measure in $\R^n$.
The quantity $\HH^n_\infty(E)$ is called $n$-dimensional Hausdorff content of $E$.
Very often we will denote the $(d-1)$-Hausdorff measure on the boundary of an open set by $\sigma$.

Recall that the fundamental solution of the minus Laplacian in $\R^d$ equals
\begin{equation*}
\EE(x)
=\begin{cases}
\dfrac{|x|^{2-d}}{(d-2)\kappa_{d}} & \mbox{if } d\geq 3, \\ &\\
\dfrac{-\log |x|}{2\pi}  & \mbox{if } d=2,
\end{cases}
\end{equation*}
for some absolute constant $\kappa_d>0$.
For a Radon measure $\mu$ in $\R^d$, we consider the energy integral
\begin{equation}\label{eq:energiamu}
I(\mu):=\iint \EE(x-y) d\mu(y)d\mu(x).
\end{equation}
Given a bounded set $E\subset\R^d$, 
we define its capacity $\capp(E)$ by
\begin{equation}\label{eqicapimu}
	\capp(E) = \frac 1{\inf_{\mu\in M_1(E)} I(\mu)},
\end{equation}
where the infimum is taken over the collection $M_1(E)$ of all {\em probability} Radon measures $\mu$ supported on $E$. When $d\geq 3$, $\capp(E)$ is also called the Newtonian capacity of $E$, and for $d=2$, the Wiener capacity of~$E$. 
Given a set $E\subset\R^2$, 
we define its logarithmic capacity by
$$\capp_L(E) = e^{-2\pi\inf_{\mu\in M_1(E)} I(\mu)}.$$

Next we recall what it means that a metric space supports Poincar\'e inequality.
 Denote by $\Sigma$  a metric space equipped with a doubling
measure $\sigma$ on $\Sigma$. That is, $\supp\sigma=\Sigma$ and there exists some constant $C>0$ such that $\sigma(B(x,2r))\leq C\,\sigma(B(x,r))$ for all $x\in\Sigma$ and $r>0$.
Given an interval $I\subset\R$, a continuous $\gamma:I \to \Sigma$ is called path. Abusing notation, we identify $\gamma$ with $\gamma(I)$.
A path of finite length is called  rectifiable path.  For a function $f: \Sigma \to \R$, we say that a Borel measurable function $g: \Sigma \to \R$ is an upper gradient of $f$ if for all compact rectifiable paths $\gamma$ the following inequality holds:
\begin{equation}\label{eq:upper-grad}
|f(x)-f(y)| \leq \int_{\gamma} g \,d\HH^1,
\end{equation}
where $x, y\in\Sigma$ are the endpoints of the path.

For $p\geq1$, we say that $\Sigma$ (or $(\Sigma,\sigma)$) supports a  weak $p$-Poincar\'e inequality if there exist constants $C\geq 1$ and $\Lambda \geq 1$ such that for every ball $B$ centered at $\Sigma$ with radius $\rad(B) \in (0, \diam (\Sigma))$ and  every  pair $(f,g)$, where $f \in L^1_{loc}(\sigma)$ and $g$ is an upper gradient of $f$, it holds that
\begin{equation}\label{eq:Poincare}
 \avint_{B} \left| f(x) - m_{\sigma,B}(f)\right| \,d\sigma(x) \leq C\, \rad(B) \left( \avint_{\Lambda B} |g(x)|^p\,d\sigma(x) \right)^{1/p},
\end{equation}
where we denoted 
$m_{\sigma,B}(f)= \,\avint_B f(y) \,d\sigma(y) = \frac1{\sigma(B)} \int_B f(y) \,d\sigma(y)$.


\section{Quantitative rectifiability and square functions.}\label{sec3}

\subsection{Jones' traveling salesman theorem.}
For applications to some questions in harmonic analysis, such as the $L^2$ boundedness of the Cauchy and Riesz transforms and the study of the removability of singularities for bounded analytic functions, it is necessary to introduce quantitative techniques for the 
study of rectifiability. This was the motivation for Peter Jones to prove what now is known as the analyst traveling salesman theorem.
To state this result we need to introduce some notation and terminology. 

For a set $E\subset\R^d$ and a cube $Q\subset\R^d$,
 we denote 
 $$\beta_{\infty,E}^n(Q) = \inf_L\, \biggl\{ \sup_{y\in E\cap Q}
\frac{\dist(y,L)}{\ell(Q)}\biggr\},
$$
where the infimum is taken over all the $n$-planes $L\subset\R^d$ and $\ell(Q)$ stands for the side length of $Q$.
Also, for a  constant $\lambda>0$, we let $\lambda Q$ be the cube concentric with $Q$ with side length $\lambda\ell(Q)$.
We denote by $\DD(\R^d)$ the family of dyadic cubes in $\R^d$.
Notice that $\beta_{\infty,E}^n(Q)$ is scale invariant and $\beta_{\infty,E}(Q)=0$ if $E\cap Q$ is contained in some $n$-plane $L$.

The Jones' traveling salesman theorem reads as follows.

\begin{theorem} \label{teojones}
A set $E\subset\R^d$ is contained in a curve of finite length if and only if
\begin{equation}\label{betajones}
\sum_{Q\in \DD(\R^d)} \beta_{\infty,E}^1(3Q)^2\ell(Q) <\infty.
\end{equation}
The length of the shortest curve $\Gamma$ containing $E$ satisfies
\begin{equation} \label{betajones2}
C^{-1}\HH^1(\Gamma) \leq \diam(E) + \sum_{Q\in \DD(\R^d)} \beta_{\infty,E}^1(3Q)^2\ell(Q) \leq C\,\HH^1(\Gamma),
\end{equation}
where $C>0$ is an absolute constant depending only on $d$.
\end{theorem}

Remark that we identify the arc-length measure with the Hausdorff measure $\HH^1$. This result was proved by Jones \cite{JonesTraveling} for subsets of the plane, while one of the implications in higher dimensions, namely that if $E$ is contained in a curve $\Gamma$ then \eqref{betajones2} holds, is due to Okikiolu \cite{Okikiolu}.

\subsection{Uniform rectifiability.}
Jones' theorem was very influential and it can be considered as the starting of the theory of quantitative rectifiability. This 
 has been extended to higher dimensions in different ways.
In particular, David and Semmes \cite{DavidSemmes1}, \cite{DavidSemmes2}, in the 1990's introduced the notion of uniform rectifiability and they characterized uniformly rectifiable
sets and measures in terms of $\beta$-type coefficients. 

Before introducing uniform rectifiability, we need to define Ahlfors regular  measures and sets.
First, we say that a Radon measure $\mu$ in $\R^d$ is $n$-Ahlfors regular (or just Ahlfors regular when $n$ is implicit) 
if there exists 
some constant $C_0>0$ such that
$$C_0^{-1}\,r^n \leq \mu(B(x,r)) \leq C_0\,r^n\quad \mbox{ for all $x\in\supp\mu$ and $0<r\leq \diam(\supp(\mu))$.}$$
If we only require the upper bound on the right hand side for $\mu$, then we say that $\mu$ has $n$-polynomial growth.
We say that $\mu$ is uniformly $n$-rectifiable if it is $n$-Ahlfors regular and there exist $\theta,M>0$ such that, for all $x\in\supp\mu$
and $0<r\leq \diam(\supp\mu)$, there exists a Lipschitz mapping from the ball $B_n(0,r)\subset\R^n$ to $\R^d$ with
$\Lip(g)\leq M$ such that
$$\mu(B(x,r)\cap g(B_n(0,r)))\geq \theta\,r^n.$$
Also, we say that a closed set $E\subset\R^d$ is $n$-Ahlfors regular if $\HH^n|_E$ is Ahlfors regular, and $E$ is uniformly $n$-rectifiable 
 if $\HH^n|_E$ is uniformly $n$-rectifiable.

For a Radon measure $\mu$ in $\R^d$, $1\leq p<\infty$, $x\in\R^d$, and $r>0$, we denote
$$\beta_{p,\mu}^n(x,r) = \inf_L \left(\frac1{r^n} \int_{B(x,r)} \left(\frac{\dist(y,L)}r \right)^p\,d\mu(y)\right)^{1/p},$$
where the infimum is taken over all $n$-planes $L\subset \R^d$. 
For $p=\infty$, we set
$$\beta_{\infty,\mu}^n(x,r) = \inf_L\, \bigg( \sup_{y\in \supp\mu}
\frac{\dist(y,L)}{r}\bigg).
$$
For a set $E\subset\R^d$, we write $\beta_{p,E}^n(x,r) = \beta_{\infty,\HH^n|_E}^n(x,r)$.

In \cite{DavidSemmes1}, David and Semmes proved the following result.

\begin{theorem}\label{teoDS}
Let $\mu$ be an $n$-Ahlfors regular measure in $\R^d$ and let $1\leq p <2n/(n-2) $.
Then, $\mu$ is uniformly $n$-rectifiable if and only if for all $x\in\supp\mu$, $0<r\leq\diam(\supp\mu)$, it holds that
$$\int_{B(x,r)}\int_0^r \beta_{\mu,p}^n(y,t)^2\,\frac{dt}t\,d\mu(y) \leq C\,\mu(B(x,r)),$$
for some constant $C>0$.
\end{theorem}

Observe that the equivalence in the theorem holds for $1\leq p\leq 2$ for any $n$. In \cite{DavidSemmes1} and \cite{DavidSemmes2}
one can find many other characterizations of uniform rectifiability. In particular, it is proven that $\mu$ is uniformly rectifiable if 
and only if any singular integral operators with an odd $n$-dimensional Calder\'on-Zygmund kernel (smooth enough) is bounded in $L^2(\mu)$.
See Section \ref{sec3} for more details. So, in a sense, uniformly $n$-rectifiable measures are the good ones to study the $L^2$ boundedness
of singular integral operators (at least, under the assumption of $n$-Ahlfors regularity). See also \cite{Tolsa-alphas} for a related result in terms of some transportation type coefficients.

\subsection{The non-Ahlfors regular case.}

It is natural to wonder if some of the results about uniform rectifiability, such as Theorem \ref{teoDS}, can be extended to non-Ahlfors regular sets or measures. The arguments in \cite{DavidSemmes1} and \cite{DavidSemmes2} depend strongly on the Ahlfors regularity condition, and so
the question is far from trivial. 
A first answer is provided by the following result.

\begin{theorem}
Let $E\subset\R^d$ be $\HH^n$-measurable and such that $\HH^n(E)<\infty$.
Then, $E$ is $n$-rectifiable if and only if  
\begin{equation}\label{eqcond49}
\int_0^1 \beta_{2,E}^n(x,r)^2\,\frac{dr}r <\infty\qquad \mbox{for $\HH^n$-a.e.\ $x\in E$.}
\end{equation}
\end{theorem}

The fact that $n$-rectifiable sets satisfy  \rf{eqcond49} was proved in \cite{Tolsa-cvpde}, and the converse implication by Azzam and the author in
\cite{AzzamTolsa}. On the other hand, as shown in \cite{Tolsa-betas}, in the theorem above one cannot replace the coefficients $\beta_{2.E}^n$ by $\beta_{p.E}^n$ for any $p\neq2$.

Later,  Naber and Valtorta \cite{NaberValtorta-Rectifiable} proved a remarkable discrete quantitative version of the above theorem which moreover provides useful estimates for
the Hausdorff measure and Minkowski content of the set $E=\supp\mu$. More interestingly, they applied their results to obtain new bounds for the size of the singular set of stationary and minimizing harmonic maps. Their result and techniques have been applied later on to study the singular sets of many other free boundary problems.  
We state below a variant of their theorem from their joint work with Edelen \cite{EdelenNaberValtorta}.

\begin{theorem}
Let $E\subset B(0,1)\subset \R^d$ be $\HH^n$-measurable and such that $\HH^n(E)<\infty$ satisfying one of the following conditions:
\begin{itemize}
\item[(a)] For $\HH^n$-a.e.\ $x\in E$,
$$\int_0^2 \beta_{E,2}^n(x,r)^2\,\frac{dr}r \leq M.$$

\item[(b)] 
For $\HH^n$-a.e.\ $x\in E$ and every $0<r\leq 1$,
$$\int_{B(x,r)\cap E} \int_0^{2r}\beta_{E,2}^n(y,t)^2\,\frac{dt}t\,d\HH^n(y) \leq M^2\,r^n.$$

\item[(c)]
For $\HH^n$-a.e.\ $x\in E$ and every $0<r\leq 1$,
$$\int_{B(x,r)\cap E} \int_0^{2r}\beta_{E,2}^n(y,t)^2\,\frac{dt}t\,d\HH^n(y) \leq M\,\HH^n(E\cap B(x,r)).$$
\end{itemize}
Then $E$ is $n$-rectifiable and for every $x\in E$, $0<r\leq1$, we have
$$\HH^n(E\cap B(x,r))\leq c(n)(1+M)r^n.$$
\end{theorem}

For other related works with a similar flavor, see David and Toro \cite{DavidToro}, Azzam and Schul \cite{AzzamSchul}, and Mi\'skiewicz \cite{Miskiewicz}. For results involving a different notion of rectifiability for measures, see for example Badger and Schul \cite{BadgerSchul}.


\subsection{Existence of tangents and Carleson $\ve^2$-conjecture.}

Next, we will define the notion of tangent in codimension $1$. To this end, first we need to introduce the notation for cones. 
For a point $x\in\R^{d}$, a unit vector $v$,
and an aperture parameter $a\in(0,1)$ we denote
$$X_a(x,v)=\bigl\{y\in\R^{d}:|(y-x)\cdot v|> a|y-x|\bigr\}.$$
Notice that this is a two sided cone with axis in the direction of $v$. 

Usually, for a set $E\subset\R^d$ and $x\in E$, one says that $x$ is a tangent point if
there exists a unit vector $v$ such that, for all
$a\in(0,1)$, there exists some $r>0$ such that
$$E\cap X_a(x,v)  \cap B(x,r) =\varnothing.$$
We need to introduce a variant of this notion for pairs of disjoint open sets: given
 $\Omega_1,\Omega_2\subset\R^{d}$ open and disjoint and $x\in\R^d$,
we say that $x$ is a tangent point for the pair $\Omega_1,\Omega_2$ if $x\in\pom_1\cap\pom_2$ and there exists a unit vector $v$ such that, for all
$a\in(0,1)$, there exists some $r>0$ such that
$$(\partial \Omega_1\cup \partial \Omega_2)\cap X_a(x,v)  \cap B(x,r) =\varnothing,$$
and moreover, one component of $X_a(x,v)\cap B(x,r)$ is contained in $\Omega_1$ and the other in $\Omega_2$.
The hyperplane $L$ orthogonal to $v$ through $x$ is called a tangent hyperplane at $x$.  In case that $\Omega_2=\R^{d}\setminus\overline{\Omega_1}$, we say that $x$ is a tangent point for $\Omega_1$.

The next result from Bishop and Jones \cite{BishopJones} from 1994 characterizes the tangents point for Jordan domains in terms of the 
$\beta_\infty^1$ coefficients. 

\begin{theorem}\label{teoBJ}
Let $\Gamma\subset\R^2$ be a Jordan curve and let $\Omega\subset\R^2$ be a Jordan domain with boundary $\Gamma$.
  Then, up to a set of zero measure $\HH^1$,
\begin{equation}\label{eqbet00}
\mbox{$\displaystyle\int_0^1 \beta_{\Gamma,\infty}^1(x,r)^2\,\frac{dr}r<\infty$\; for $x\in\Gamma$ \quad$\Leftrightarrow$ \quad
$\Omega$ has a tangent at $x\in\Gamma$.}
\end{equation}
  
  \end{theorem}

The proof of this theorem in \cite{BishopJones} relies strongly on Theorem \ref{teojones}.

The study of a two-phase problem for harmonic measure (see also Section \ref{sec6}), led Bishop, Carleson, Garnett, and Jones to
introduce another square functions connected both to rectifiability and to harmonic measure. This is the so-called Carleson's $\ve^2$-square function. To define this, let $\Omega_1$ be a Jordan domain in $\R^2$, and set $\Gamma =\partial\Omega_1$ and
$\Omega_2 = \R^2\setminus \overline{\Omega_1}$.
For $x\in\R^2$ and $r>0$, denote by $I_1(x,r)$ and $I_2(x,r)$ the longest open arcs of the circumference
$\partial B(x,r)$ contained in $\Omega_1$ and $\Omega_2$, respectively (they may be empty). Then we set
\begin{equation}\label{eqepsiloncoef}
\ve(x,r) = \frac1r\,\max\big(\big|\pi r- \HH^1(I_1(x,r))\big|,\, \big|\pi r- \HH^1(I_2(x,r))\big|\big).
\end{equation}
In particular, observe that $\ve(x,r)$ vanishes if and only if the arcs $I_1(x,r)$ and $I_2(x,r)$ are semicircumferences (i.e., they correspond to plane angles).
  The Carleson $\ve^2$-square function is defined by
\begin{equation}\label{eqEE}
\EE_C(x)^2 :=\int_0^1 \ve(x,r)^2\,\frac{dr}r.
\end{equation}
One may expect that this function encodes some information related to rectifiability.
In fact,  Carleson's conjecture, now a theorem, asserts the following.

\begin{theorem} \label{teo-carleson}
Let $\Omega_1\subset \R^2$ be a Jordan domain, let $\Omega_2=\R^2\setminus \overline{\Omega_1}$, and let $\Gamma=
\partial\Omega_1$. Let $\EE_C$ be the
associated square function defined in \rf{eqEE}. Then, up to  set of zero measure $\HH^1$,
\begin{equation}\label{eqCar00}
\mbox{$\displaystyle\EE_C(x)^2 =\int_0^1 \ve(x,r)^2\,\frac{dr}r<\infty$\; for $x\in\Gamma$ \quad$\Leftrightarrow$ \quad
$\Omega$ has a tangent at $x\in\Gamma$.}
\end{equation}
In particular, the set
$G=\{x\in\Gamma:\EE(x)<\infty\}$ is $1$-rectifiable.
\end{theorem}

The fact  that $\EE_C(x)<\infty$ for $\HH^1$-a.e.\ tangent point in a Jordan curve was proved by Bishop in 
\cite{Bishop-Thesis} (see also \cite{BCGJ}). The most difficult implication of Theorem \ref{teo-carleson}, i.e, the fact that the set $G$ is $1$-rectifiable and tangents to $\Gamma$ exist for $\HH^1$-a.e.\ $x\in G$,
was proved more recently by Jaye, the author, and Villa \cite{JayeTolsaVilla}.

The arguments from \cite{Bishop-Thesis} and \cite{BCGJ} which show the finiteness of the square function $\EE_C(x)$ at tangent points in $\Gamma$ use harmonic measure estimates. To explain them, denote by $\omega_i^{p_i}$ the harmonic measure for $\Omega_i$ with pole at $p_i\in\Omega_i$ (see Section \ref{sec5} for more details).
  Already in the 1930's Ahlfors obtained the following  distortion result for harmonic measure in Jordan domains: \begin{equation}\label{beaurling0}
	\omega_i^{p_i}(B(x, r)) \leq C \exp \left(-\pi \int_r^R \frac{1}{\Theta_i(x,t)} \, \frac{dt}{t}\right),\quad
	\mbox{ for  $p_i\in\Omega_i$, $x \in \partial \Omega_i$, and $0<r<R\leq |p_i-x|$,}
\end{equation}
where $C$ is some constant depending on $\dist(p,\pom_i)$ and $R$, and, for $I_i(x, t)$ being the longest open arc in $\Omega_i \cap \partial B(x,t)$,  $\Theta_i(x,t)= \HH^1(I_i(x,t))/t$.\footnote{See Chapter 6 from Garnett and Marshall \cite{GarnettMarshall}.}

Under the assumptions in the theorem,  we may use \eqref{beaurling0} to see that
\begin{equation}\label{beaurling1}
	\frac{\omega_1^{p_1}(B(x, r))\,\omega_2^{p_2}(B(x,r))}{r^2} \lesssim \exp \left(- \pi\int_{r}^R \left(\frac{1}{\Theta_1(x,t)} + \frac{1}{\Theta_2(x,t)}-\frac2\pi\right) \, \frac{dt}{t} \right).
\end{equation}
By McMillan's theorem, the harmonic measures $\omega_i^{p_i}$ (for $p_i \in \Omega_i$) are mutually absolutely continuous with the Hausdorff measure $\HH^1$ on the subset of tangent points (see Theorem \ref{teoconepoints} below), and since this set is rectifiable, we deduce that the densities 
$$\theta(x,\omega_i^{p_i}) := \lim_{r\to0}\frac{\omega_i^{p_i}(B(x,r))}{2r}$$
exist and are positive and finite, for $i=1,2$, at $\HH^1$-a.e.\ tangent point $x$.
 So letting $r\to0$ in \eqref{beaurling1}, we deduce that
\begin{equation}\label{finite1}
	\int_{0}^R \left(\frac{1}{\Theta_1(x,t)} + \frac{1}{\Theta_2(x,t)} -\frac2\pi\right) \, \frac{dt}{t} < + \infty
\end{equation}
for $\HH^1$-a.e.\ tangent point $x\in\pom_i$.
Now, a Taylor expansion shows that 
\begin{equation*}
	\frac{1}{\Theta_1(x,t)} + \frac{1}{\Theta_2(x,t)} \geq \frac{2}{\pi} + \frac{2}{\pi} \left( \frac{\ve(x, t)}{\pi}\right)^2,
\end{equation*}
where $\ve(x,t)$ is as in \rf{eqepsiloncoef}.
In view of \eqref{beaurling1}, we see that at $\HH^1$-almost all tangent points of $\Gamma$, $\EE_C(x)<\infty$.

An alternative argument to prove the finiteness of the square function $\EE_C(x)$ at tangent points in $\Gamma$ in Theorem \ref{teo-carleson} is based 
on the application of Theorem \ref{teoBJ}.
Indeed, it is easy to check that  whenever $x$ is a tangent point of $\Omega_1$ and r is small enough, it holds that
$\ve(x,r) \leq c\,\beta_{\Gamma,\infty}^1(x,r)$. So one deduces that $\EE_C(x)<\infty$ for $\HH^1$-a.e.\ tangent point in $\Gamma$ as an immediate corollary of Theorem \ref{teoBJ}.\footnote{Remark that the work \cite{BishopJones} where Theorem \ref{teoBJ} is proved is from 1994 while \cite{Bishop-Thesis} and
\cite{BCGJ} are from 1987.} However, the converse inequality 
between $\ve(x,r)$ and $\beta_{\Gamma,\infty}^1(x,r)$ fails
and thus one cannot 
deduce the other implication in Theorem \ref{teo-carleson} from Theorem \ref{teoBJ}. Instead, the arguments in \cite{JayeTolsaVilla} combine some  delicate blowup arguments with a stopping time construction inspired in part by the techniques from \cite{DavidSemmes1} and \cite{Leger}.

\subsection{Carleson $\ve^2$-conjecture in higher dimensions.}

Recently, in \cite{FleschlerTolsaVilla-Inventiones}, Fleschler, the author, and Villa,  have obtained a higher dimensional 
 version of Theorem \ref{teo-carleson}. To state this, we need some additional terminology.
 For $x \in \R^{d}$, $r>0$ and an open half-space $H$ such that $x \in \partial H$, we denote 
\begin{equation}\label{eqSH}
S_{H}^1(x,r)=S(x,r) \cap H,\qquad S_{H}^2(x,r):= S(x,r) \cap (\R^{d} \setminus \overline{H}).
\end{equation}
Given two disjoint Borel sets $\Omega_1,\Omega_2\subset\R^{d}$, put
\begin{equation}\label{ve_n-intro}
	\ve_{d-1}(x,r) := \frac{1}{r^{d-1}}\, \inf_H \HH^{d-1} \left( (S_H^1(x,r) \setminus \Omega_1) \cup (S_H^2(x,r) \setminus \Omega_2)\right),
\end{equation}
where the infimum is taken over all half-spaces $H$ such that $x\in\partial H$.
It is clear that if $\Omega_1$ and $\Omega_2$ are complementary (open) half-spaces, then $\ve_{d-1}(x,r)=0$ for any  $x \in \partial \Omega_1 = \partial \Omega_2$ and $r>0$. Note that in the plane $\ve_1(x,r)  \lesssim \ve(x,r)$, but the opposite inequality fails, in general. 
We write 
\begin{equation}\label{calE_n}
	\mathcal{E}_{d-1}(x)^2 := \int_0^1 \ve_{d-1}(x,r)^2 \, \frac{dr}{r}.
\end{equation}
The first main result from \cite{FleschlerTolsaVilla-Inventiones} is the following:

\begin{theorem}\label{teocarleson2}
	For $d\geq 2$, let $\Omega_1, \Omega_2 \subset \R^{d}$ be two disjoint Borel subsets. Then any Borel subset of $\{x \in \R^{d}: \mathcal{E}_{d-1}(x) < \infty\}$ is ${d-1}$-rectifiable. 
\end{theorem}

Remark that the assumptions in this theorem are much weaker than the ones in Theorem \ref{teo-carleson} in the planar case.
  However, Theorem \ref{teocarleson2} does not ensure the existence of tangents for the pair of sets $\Omega_1$, $\Omega_2$. 
As in \cite{JayeTolsaVilla}, the proof combines blowup arguments with quantitative rectifiability techniques, although they are
 more elaborate and difficult than the ones from the planar case.

To ensure the existence of tangents for $\Omega_1$, $\Omega_2$ one needs stronger assumptions and different coefficients from the $\ve_{d-1}$'s. These new coefficients  involve the Dirichlet eigenvalues of some open sets in the unit sphere.
Given an open set $V$  from $\bS^{d-1}$, one says that $u\in W^{1,2}_0(V)$ is a Dirichlet eigenfunction of
$V$ for the Laplace-Beltrami operator $-\Delta_{\bS^{d-1}}$ if $u\not\equiv 0$ and
$$-\Delta_{\bS^{d-1}} u = \lambda\,u,$$
for some $\lambda\in \R\setminus \{0\}$. Here $W^{1,2}_0(V)$ is the closure in the Sobolev space $W^{1,2}(V)$ of the $C^\infty$ functions  with compact support in $V$.
The number $\lambda$ is the eigenvalue associated with $u$. It is well known that all the eigenvalues of the Laplace-Beltrami operator are positive.
The characteristic constant of $V$ is defined as the positive number $\alpha_V$ such that $\lambda_V = \alpha_V(d-2+\alpha_V)$, where $\lambda_V$ is the smallest eigenvalue.

The Friedland-Hayman \cite{FriedlandHayman} inequality asserts that, for any two disjoint open subsets $V_1,V_2\subset \bS^{d-1}$ and $\alpha_i\equiv
\alpha_{V_i}$, it holds that
 $$\alpha_1+\alpha_2-2\geq 0.$$
 This inequality plays an essential role in the proof of 
the Alt-Caffarelli-Friedman (ACF) monotonicity formula, which is a basic tool in many free boundary problems. 
We recall that this asserts the following.

	\begin{theorem} \label{teoACF-elliptic}  Let $x \in  \R^{d}$ and $R>0$. Let $u_1,u_2\in
		W^{1,2}(B(x,R))\cap C(B(x,R))$ be nonnegative subharmonic functions such that $u_1(x)=u_2(x)=0$ and $u_1\cdot u_2\equiv 0$. 
		Set
		\begin{equation*}
			J(x,r) = \left(\frac{1}{r^{2}} \int_{B(x,r)} \frac{|\nabla u_1(y)|^{2}}{|y-x|^{d-2}}dy\right)\cdot \left(\frac{1}{r^{2}} \int_{B(x,r)} \frac{|\nabla u_2(y)|^{2}}{|y-x|^{d-2}}dy\right)
		\end{equation*}
		Then $J(x,r)$ is an absolutely continuous function of $r\in (0,R)$ and
			\begin{equation}\label{eqprec1}
			\frac{\partial_rJ(x,r)}{J(x,r)}\geq \frac2r\bigl(\alpha_1 + \alpha_2  - 2 \bigr). 
		\end{equation}
		where $\alpha_i$ is the characteristic constant of the open subset  $V_i\subset\bS^{d-1}$ given by  $$V_i=\bigl\{r^{-1}(y-x): y\in\partial B(x,r),\,u_i(y)>0\bigr\}.$$
	\end{theorem}
	
Notice that the Friedland-Hayman inequality together with \rf{eqprec1} ensure that $J(x,r)$ is non-decreasing in $r$, under the assumptions of the 
theorem.

Given two disjoint open sets $\Omega_1,\Omega_2\subset \R^{d}$ and $x\in\R^{d}$, $r>0$, we consider
the sets
$V_i(x,r) =\{r^{-1}(x-y) \, : \, y \in S(x,r)\cap \Omega_i\}$ and we denote
\begin{equation}\label{eqalfai}
\alpha_i(x,r) = \alpha_{V_i(x,r)}.
\end{equation}
Notice that $\alpha_1(x,r) + \alpha_2(x,r) - 2\geq0$ by the Friedland-Hayman  inequality.
In the plane, an easy computation (see \cite{AllenKriventsovNeumayer-ACF}) shows that the coefficients $\ve(x,r)$ defined in \rf{eqepsiloncoef} satisfy
\begin{equation}\label{eqalfai2}
\ve(x,r)^2\lesssim \min\big(1,\alpha_1(x,r) + \alpha_2(x,r) - 2\big) .
\end{equation}
Further, in case that $I_1(x,r)$ and $I_2(x,r)$ are complementary arcs in $S(x,r)$, one can check that
$$\min\big(1,\alpha_1(x,r) + \alpha_2(x,r) - 2\big) 
\approx \ve(x,r)^2.$$

In higher dimensions we have 
\begin{equation}\label{eqclau10}
\ve_{d-1}(x,r)^2 \lesssim \min(1,\alpha_1(x,r) + \alpha_2(x,r)-2).
\end{equation}
This was proved in \cite{FleschlerTolsaVilla-Faber}. The arguments are much more delicate than the ones for \rf{eqalfai2} and they involve suitable quantitative Faber-Krahn inequalities arising from remarkable works of Brasco, De Philippis, and Velichkov \cite{BrascoPhilippisVelichkov} and 
Allen, Kriventsov, and Neumayer \cite{AllenKriventsovNeumayer-Ars}. 

In \cite{FleschlerTolsaVilla-Inventiones} we have obtained a version of Carleson's conjecture in terms of the square function
\begin{equation*}
\mathcal A(x)^2 := \int_0^1 \min(1,\alpha_1(x,r) + \alpha_2(x,r)-2)\,\frac{dr}r.
\end{equation*} 
The precise result is the following.

\begin{theorem}\label{teomain2}
	For $d\geq 2$ let $\Omega_1, \Omega_2 \subset \R^{d}$ be two disjoint open subsets. Suppose that $\Omega_1 \cup \Omega_2$ satisfies the capacity density condition (CDC). Then, up to a set of zero $\HH^{d-1}$ measure,
	\begin{equation*}
		\mathcal{A}(x) < \infty \,\, \mbox{ if and only if $x$ is a tangent point of the pair} \,\, \Omega_1, \Omega_2.
	\end{equation*} 
\end{theorem}

The CDC for an open set $\Omega$ (such as $\Omega_1 \cup \Omega_2$) is a well known condition in potential theory which requires some thickness
for $\Omega^c$. More precisely, a set $\Omega\subset \R^d$, with $d\geq3$ satisfies the CDC if there exists some $c>0$ such that for all $\xi\in\pom$ and $0<r\leq\diam(\pom)$,
$$\capp(B(\xi,r)\setminus \Omega)\geq c\,r^{d-2},$$
where $\capp$ is the Newtonian capacity of homogeneity $d-2$. In the case $d=2$, instead one asks that for all $\xi\in\pom$ and $0<r\leq\diam(\pom)$,
$$\capp_L(B(\xi,r)\setminus \Omega)\geq c\,r,$$
where $\capp_L$ is the logarithmic capacity.
For example, if $\Omega_1\subset\R^2$ is a Jordan domain and $\Omega_2=\R^2\setminus \overline{\Omega_1}$, then the CDC holds $\Omega_1 \cup \Omega_2$. So Theorem \ref{teomain2} implies Theorem \ref{teo-carleson}, and in fact it extends to more general open sets in the plane.

The fact that the square function $\mathcal{A}$ is finite at tangent points of the pair $\Omega_1,\Omega_2$ follows by  a rather quick application of the ACF monotonicity formula. The converse implication is more difficult. The proof uses Theorem \ref{teocarleson2} and a version of \rf{eqclau10} for coefficients $\ve_s(x,r)$ of higher codimension from \cite{FleschlerTolsaVilla-Faber}.

In the case when $\partial\Omega_1\cup \partial\Omega_2$ is $(d-1)$-Ahlfors regular, Casey, the author, and Villa \cite{CaseyTolsaVilla} have obtained a quantitative version of the results above in terms of the so-called corkscrew condition. An open set $\Omega$ is said to satisfy the corkscrew condition if for every ball $B$ centered at $\Omega$ with radius $0<\rad(B)\leq \diam(\Omega)$, there exists a ball $B'\subset B\cap\Omega$ such that $\rad(B')\geq c \,\rad(B)$, for some fixed $c>0$. The result reads as follows.

\begin{theorem}
	Let $\Omega_1$ and $\Omega_2$ be two disjoint open subsets of $\R^{d}$. Suppose that $\mu$ is a $(d-1)$-Ahlfors measure with $\supp(\mu) = \partial \Omega_1 \cup \partial \Omega_2$. Then the following are equivalent.
	\begin{enumerate}
		\item[(a)] $\Omega_1$ and $\Omega_2$ are complementary corkscrew open sets, and in particular $\mu$ is uniformly $(d-1)$-rectifiable.
		\item[(b)] There is a constant $C$ such that for each ball $B$ centered on $\supp(\mu)$ with $\rad(B)\leq \diam(\pom_1)$ it holds that
		\begin{equation*}
			\int_B \int_0^{\rad(B)} \ve_{d-1} (x,r)^2 \, \tfrac{dr}{r} \, d\mu(x) \leq C\, \rad(B)^{d-1}.
		\end{equation*}
		\item[(c)] There is a constant $C$ such that for each ball $B$ centered on $\supp(\mu)$ with $\rad(B)\leq \diam(\pom_1)$ it holds that
		\begin{equation*}
			\int_B \int_{0}^{\rad(B)} \min(1,\alpha_1(x,r) + \alpha_2(x,r)-2) \, \tfrac{dr}{r} \, d\mu(x) \leq C \,\rad(B)^{d-1}.
		\end{equation*}
	\end{enumerate}
\end{theorem}

Regarding the condition (a) above, we remark that if $\Omega\subset\R^d$ is an open set with $(d-1)$-Ahlfors regular boundary and both $\Omega$ and $\R^d\setminus \overline{\Omega}$ 
are corkscrew domains, then $\pom$ is $(d-1)$-uniformly rectifiable. This was proved independently by David and Jerison \cite{DavidJerison} and Semmes \cite{Semmes} while studying some questions on harmonic measure. See Section \ref{sec5}.
\vv


\section{Riesz transforms and rectifiability.} \label{sec4}

Let $0<n\leq d$ and let $\nu$ be a signed Radon measure in $\R^d$. The $n$-dimensional Riesz transform of $\nu$ (or $n$-Riesz transform) 
is defined by 
$$\RR^n\nu(x) =\int\frac{x-y}{|x-y|^{n+1}}\,d\nu(y),$$
for any $x\in\R^d$ where the integral makes sense.
For any $\ve>0$, the $\ve$-truncated version of the $n$-Riesz transform is defined by 
$$\RR_\ve^n\nu(x) =\int_{|x-y|>\ve}\frac{x-y}{|x-y|^{n+1}}\,d\nu(y).$$
We also consider the principal value $n$-Riesz transform and the maximal $n$-Riesz transform, defined respectively by
$$\pv \RR^n\nu(x) = \lim_{\ve\to0} R_\ve^n\nu(x),\qquad \RR_*^n\nu(x) = \sup_{\ve >0} |R_\ve^n\nu(x)|.$$
For a positive Radon measure $\mu$ and a function $f\in L^1_{loc}(\mu)$, we write
$$\RR_\mu^n f = \RR^n(f\mu),\quad \RR_{\mu,\ve}^n f = \RR_\ve^n(f\mu),\quad \pv \RR_\mu^n f = \pv \RR^n(f\mu),\quad \RR_{\mu,*}^n f = \RR_*^n(f\mu).$$
We say that the $n$-Riesz transform is bounded in $L^2(\mu)$ if the operators $\RR_{\mu,\ve}^n$ are bounded uniformly on $\ve>0$.
This definition avoids the issue of the existence of $\RR_\mu^n f(x)$ or $\pv \RR_\mu^n f(x)$, which is not guaranteed.

One may also consider more general operators of the form
$$T\nu(x) = \int K(x,y)\,d\nu(y),$$
where $K:\R^d\times\R^d\setminus \{x=y\}\to\R$ is a Calder\'on-Zygmund kernel satisfying
$$|K(x,y)|\leq C\,\frac1{|x-y|^n}\quad \mbox{ for all $x,y\in\R^d$ with $x\neq y$}$$
and
$$|K(x,y)-K(x',y)| + |K(y,x)-K(y,x')| \leq C\,\frac{|x-x'|^\delta}{|x-y|^{n+\delta}}\quad \mbox{ for all $x,x',y\in\R^d$ such that $|x-x'|\leq \frac12|x-y|,$}$$
where $C>0$ and $0<\delta\leq 1$ are some fixed constants. Abusing notation, we also write $K(x,y) = K(x-y)$ in case $K$ is a convolution kernel.
The Riesz kernel $x/|x|^{n+1}$ is a vectorial example of this
type of kernels. One defines $T_\ve \nu$, $T_\mu f$, $T_{\mu,\ve} f$, $\pv T$ , and $T_*$ in an analogous way as in the case of the Riesz transform. This kind of operators
are often called singular integral operators.

One of the motivations of David and Semmes to develop the theory of uniform rectifiability was the wish to find a natural class of sets and measures $\mu$ such that the Riesz transform $\RR_\mu^n$ and other singular integral operators with an odd Calder\'on-Zygmund kernel are bounded in $L^2(\mu)$. In turn, it was expected (and this was the case) that these results might have applications to other problems in analysis, such as the Painlev\'e problem and Vitushkin's conjecture. Indeed, notice that, modulo a constant factor, the $(d-1)$-Riesz kernel $\frac{x}{|x|^d}$ coincides with gradient of the fundamental solution of the Laplacian, denoted by $\EE$.
 So $\RR^{d-1}_\mu f(x) =c\int\nabla \EE(x-y)\,f(y)\,d\mu(y)$, and then it is natural to think 
that the $(d-1)$-Riesz transform may play an important role in many topics 
connected to the Laplace equation and harmonic functions. Remark also that in $\R^2$, the kernel $z/|z|^2$ coincides with the conjugated of the Cauchy kernel $1/z$ and so $\RR^1$ coincides with the conjugated Cauchy transform.

In \cite{DavidSemmes1}, David and Semmes proved the following.

\begin{theorem}\label{teorieszDS}
Let $\mu$ be an $n$-Ahlfors regular measure in $\R^d$. Then, $\mu$ is uniformly $n$-rectifiable if and only if
any singular integral operator $T_\mu$ associated with an odd Calder\'on-Zygmund kernel $K$ being $C^\infty$ away from the origin and satisfying
$|\nabla^j K(x)|\leq C(j)\,|x|^{-n-j}$ for every $j\geq0$ is bounded in $L^2(\mu)$. 
\end{theorem}

In particular, the $n$-Riesz transform $\RR^n_\mu$ is bounded in $L^2(\mu)$ if $\mu$ is uniformly $n$-rectifiable. An obvious question is if the $L^2(\mu)$
boundedness of the $n$-Riesz transform alone suffices for the uniform $n$-rectifiability of $\mu$, assuming $\mu$ to be $n$-Ahlfors regular.
This is the so called David-Semmes problem. A related question of more qualitative nature without the $n$-Ahlfors regularity assumption is the following: for a set $E\subset\R^d$ with $\HH^n(E)<\infty$ and $\mu=\HH^n|_E$, does the $L^2(\mu)$ boundedness of $\RR^n_\mu$ imply the $n$-rectifiability of $E$?

The David-Semmes problem and its qualitative variant are solved in the case of dimension one ($n=1$) and codimension one ($n=d-1$):

\begin{theorem}\label{teoriesz1}
Let $n=1$ or $n=d-1$.
Let $\mu$ be an $n$-Ahlfors regular measure in $\R^d$ such that the $n$-Riesz transform $\RR_\mu^n$ is bounded in $L^2(\mu)$.
Then $\mu$ is uniformly $n$-rectifiable.
\end{theorem}

The case $n=1$ of the theorem was solved in 1996 by Mattila, Melnikov, and Verdera \cite{MattilaMelnikovVerdera} using the connection between Menger curvature and the Cauchy kernel. The case of codimension one for $d>2$ is due to Nazarov, the author, and Volberg \cite{NazarovTolsaVolberg1}. The proof uses the so called BAUP criterion for uniform rectifiability of David and Semmes \cite{DavidSemmes2},
some quasiorthogonality techniques, and a variational argument which has its roots on a previous result
by Eiderman, Nazarov and Volberg \cite{EidermanNazarovVolberg}, where it was shown that if a set $E\in\R^d$ with $0<\HH^{d-1}(E)<\infty$ satisfies
$$\theta_*^{d-1}(x,E):= \liminf_{r\to0}\frac{\HH^{d-1}(E\cap B(x,r))}{(2r)^{d-1}}=0\quad \mbox{ for $\HH^{d-1}$-a.e.\ $x\in E$,}$$
then the Riesz transform $\RR_{\HH^{d-1}|_E}^{d-1}$ cannot be bounded in $L^2(\HH^{d-1}|_E)$. That variational argument, together with the usual maximum
principle for harmonic functions, allows one to transfer some key $L^\infty$ bounds for $\RR^{d-1}(\HH^{d-1}|_E)$ on $E$ to the ambient space, where one can use the Fourier transform to obtain estimates from below for $\RR^{d-1}(\HH^{d-1}|_E)$.

In connection with the qualitative version of the David-Semmes problem, we have:

\begin{theorem}\label{teoriesz2}
Let $n=1$ or $n=d-1$.
Let $E\subset \R^d$ be a Borel set such that $\HH^n(E)<\infty$ and set $\mu =\HH^n|_E$. If $\RR_\mu^n$ is bounded in $L^2(\mu)$, then
$\mu$ is $n$-rectifiable.
\end{theorem}

In the case $n=1$, this result was proved in \cite{Leger} by L\'eger  and in a previous unpublished manuscript by Guy David.  The case $n=d-1$, $d>2$, is again due to Nazarov, the author, and Volberg \cite{NazarovTolsaVolberg2}.

The David-Semmes problem is still open  when $n$ is an integer such that $1<n<d-1$:

\begin{question}
For an integer $n$ such that $1<n<d-1$, if $\mu$ is $n$-Ahlfors regular and $\RR_\mu^n$ is bounded in $L^2(\mu)$, is then $\mu$ uniformly
$n$-rectifiable?
\end{question}

 The main obstacle to answer this question by the same techniques as in
\cite{NazarovTolsaVolberg1} is the absence of a suitable maximum principle for functions of the form $\RR^n\eta$ for $n<d-1$.
On the other hand, as shown in \cite{Tolsa-jfa}, for any integer $n$ with $1\leq n\leq d-1$, given 
 a Borel set $E\subset \R^d$ such that $\HH^n(E)<\infty$, the existence of the principal values $\pv \RR^n(\HH^n|_E)(x)$ for $\HH^n$-a.e.\
$x\in E$ implies the $n$-rectifiability of $E$. 
Under the additional assumption that
$\theta_*^n(x,E)>0$ for $\HH^n$-a.e.\ $x\in E$,
this had been proved previously by Mattila and Preiss \cite{MattilaPreiss}.

As a corollary of Theorem \ref{teoriesz2} one gets the following criterion for rectifiability.

\begin{corollary}\label{coromaximal}
Let $n=1$ or $n=d-1$.
Let $E\subset \R^d$ be a Borel set such that $\HH^n(E)<\infty$. Let $\mu$ be a Radon measure such that $\mu$ is mutually absolutely continuous with $\HH^n$ on $E$.
 Then $E$ is $n$-rectifiable if and only if 
$$\RR^{n}_*\mu(x)<\infty \quad \mbox{ for $\mu$-a.e.\ $x\in E$.}$$
\end{corollary}

Using some techniques related to the ones for the solution of the David-Semmes problem, in \cite{GirelaTolsa}, Girela-Sarri\'on and the author
obtained a criterion for the absolute continuity of a measure with polynomial growth with respect to Hausdorff measure.
To announce it, we use the following notation. For a ball $B$ and a measure $\mu$, we write
$$\Theta_\mu^n(B) = \frac{\mu(B)}{\rad(B)^n}$$
and
$$P_\mu^n(B) = \sum_{j\geq0} 2^{-j}\,\Theta_\mu^n(2^jB).$$
For $B(x,r)$, we write $\Theta_\mu^n(x,r) = \Theta_\mu^n(B(x,r))$ and $P_\mu^n(x,r) = P_\mu^n(B(x,r))$. For an $n$-plane $L\subset\R^d$,
we let
$$\beta_{1,\mu}^{n,L}(B) = \frac1{r^n} \int_{B} \frac{\dist(y,L)}r \,d\mu(y),$$
so that $\beta_{1,\mu}^{n,L}(B)= \inf_L\beta_{1,\mu}^{n,L}(B)$.

Instead of announcing the precise result from \cite{GirelaTolsa}, we state a slight variant derived in \cite{AzzamMourgoglouTolsa-2phase}:

\begin{theorem} \label{teoGT}
Let $\mu$ be a Radon measure in $\R^{d}$ and let $B\subset \R^{d}$ be a ball with $\mu(B)>0$ such that the following conditions
hold:
\begin{itemize}
\item[(a)] For some constant $C_0>0$, $P_\mu^{d-1}(B) \leq C_0\,\Theta_\mu^{d-1}(B)$.

\item[(b)] There is some $(d-1)$-plane $L$ passing through the center of $B$ such that, for some constant $0<\delta\ll 1$, $\beta_{\mu,1}^{d-1,L}(B)\leq \delta\,\Theta^{d-1}_\mu(B)$.

\item[(c)] For some constant $C_1>0$, there is $G_B\subset B$
such that
$$\sup_{0<r\leq 2 \rad(B)} \Theta_\mu^{d-1}(x,r) + \RR^{d-1}_*(\chi_{2B}\,\mu)(x)\leq 
C_1\,\Theta^{d-1}_\mu(B)\quad \mbox{ for all $x\in G_B$}$$
and
$$\mu(B\setminus G_B)\leq \delta \,\mu(B).$$

\item[(d)] For some constant $0<\tau\ll1$,
$$\int_{G_B} |\RR^{d-1}\mu(x) - m_{\mu,G_B}(\RR^{d-1}\mu)|^2\,d\mu(x) \leq \tau \,\Theta^{d-1}_\mu(B)^2\mu(B).$$
\end{itemize}

Then there exists some constant $\theta>0$ such that if $\delta,\tau$ are small enough (depending on $C_0$ and $C_1$),
then there is a uniformly $(d-1)$-rectifiable set $\Gamma\subset\R^{d}$ such that
$$\mu(G_B\cap \Gamma)\geq \theta\,\mu(B).$$
The uniform rectifiability constants of $\Gamma$ depend on all the constants above.
\end{theorem}

We remark that this theorem plays a key role in the solution of the two-phase problem for harmonic measure in \cite{AzzamMourgoglouTolsa-2phase} and
\cite{AzzamMourgoglouTolsaVolberg}. See Section \ref{sec6} for more details.

The condition (a) above should be understood as a kind of doubling condition for the ball $B$. On the other hand, (b) asks the measure in $B$ to be concentrated close some hyperplane containing the center of the ball. The condition (d) is a key assumption that requires a small oscillation
of $\RR^{d-1}\mu$ with respect to its mean $\mu$ on $G_B$, denoted by $m_{\mu,G_B}(\RR^{d-1}\mu)$.
 It is not clear to the author whether the assumptions in the theorem are sharp.
For example, it is natural to wonder if the assumption (b) can be eliminated.
 
By 
Theorems \ref{teorieszDS} and \ref{teoriesz1} we have a complete metric-geometric characterization of the $(d-1)$-Ahlfors regular measures such that the $(d-1)$-Riesz transform is bounded in $L^2$: they are the uniformly $(d-1)$-rectifiable measures. For applications to questions such as the Painlev\'e
problem for Lipschitz harmonic functions, it is desirable to have a characterization for arbitrary Borel measures with no point masses. This has been achieved
by D\k{a}browski and the author in \cite{DabrowskiTolsa}. The precise result is the following.

\begin{theorem}\label{teomain22}
Let $\mu$ be a Radon measure in $\R^{d}$ with no point masses. Then $\RR^{d-1}_\mu$ is bounded in $L^2(\mu)$ if and
only if it satisfies the polynomial growth condition
\begin{equation}\label{eqgrow01}
\mu(B(x,r))\leq C\,r^{d-1}\quad \mbox{ for all $x\in\supp\mu$ and all $r>0$}
\end{equation}
and
\begin{equation}\label{eqbetawolff2}
\int_B\int_0^{\rad(B)} \beta_{2,\mu}^{d-1}(x,r)^2\,\Theta_\mu^{d-1}(x,r)\,\frac{dr}r\,d\mu(x)\leq C^2\,\mu(B)\quad\mbox{ for any ball
$B\subset\R^{d}$.}
\end{equation}
Further, the optimal constant $C$ is comparable to $\|\RR^{d-1}_\mu\|_{L^2(\mu)\to L^2(\mu)}$.
\end{theorem}

\vv
Moreover, by combining the results from \cite{DabrowskiTolsa} and the ones of Girela-Sarri\'on in \cite{GirelaSarrion},  under the assumption that
$\mu$ has $(d-1)$-polynomial growth and that $\|\RR_*^{d-1}\mu\|_{L^1(\mu)}<\infty$,
it follows that
$$
\int\!\!\int_0^\infty \beta^{d-1}_{2,\mu}(x,r)^2\,\Theta^{d-1}_\mu(x,r)\,\frac{dr}r\,d\mu(x)+ C\,\|\mu\|\approx\|\pv\RR^{d-1}\mu\|_{L^2(\mu)}^2
+C\,\|\mu\|.$$
In the case $d=2$, Theorem \ref{teomain22} had been obtained previously in \cite{AzzamTolsa}, relying on the connection between the Cauchy transform and Menger curvature.

As shown in \cite{AzzamTolsa}, the condition \rf{eqbetawolff2}
is equivalent to the existence of a suitable corona decomposition for $\mu$ satisfying an appropriate packing
condition. This condition is stable by bilipschitz maps (see also \cite{GirelaSarrion} for more details). So we get the
following corollary. Again, this was already known  previously in the planar case $d=2$ \cite{Tolsa-bilip}.

\begin{corollary}\label{coro1}
Let $\mu$ be a Radon measure in $\R^{d}$ with no point masses. Let $\vphi:\R^{d}\to\R^{d}$ be a bilipschitz
map. Let $\sigma=\vphi\#\mu$ be the image measure of $\mu$ by $\vphi$.
If $\RR_\mu^{d-1}$ is bounded in $L^2(\mu)$, then $\RR_\sigma^{d-1}$ is bounded in $L^2(\sigma)$. Further,
$$\|\RR_\sigma^{d-1}\|_{L^2(\sigma)\to L^2(\sigma)}\leq C\,\|\RR_\mu^{d-1}\|_{L^2(\mu)\to L^2(\mu)},$$
where $C$ depends only on the bilipschitz constant of $\vphi$.
\end{corollary}

The image measure  $\vphi\#\mu$ is the measure defined by $\vphi\#\mu(A) = \mu(\vphi^{-1}(A))$ for any Borel set $A$. 
We remark that, up to now, the preceding result was not  known even for the case of invertible affine maps such as the one defined by
$\vphi(x_1,\ldots,x_{d}) = (2x_1,x_2,\ldots,x_{d})$.

The conditions \rf{eqgrow01} and \rf{eqbetawolff2} imply the $L^2(\mu)$ boundedness of any singular integral operator 
$T_\mu$ associated with an odd Calder\'on-Zygmund kernel $K$ satisfying
$|\nabla^j K(x)|\leq C(j)\,|x|^{-n-j}$ for  $0\leq j\leq 2$. See  \cite{GirelaSarrion}.
Then we deduce the following.

\begin{corollary}\label{coro1.5}
Let $\mu$ be a Radon measure in $\R^{d}$ with no point masses. Let $T_\mu$ be a singular integral 
operator associated with an odd kernel $K$ satisfying $|\nabla^j K(x)|\leq C(j)\,|x|^{-d+1-j}$ for  $0\leq j\leq 2$.
If $\RR_\mu^{d-1}$ is bounded in $L^2(\mu)$, then $T_\mu$ is also bounded in $L^2(\mu)$. Further,
$$\|T_\mu\|_{L^2(\mu)\to L^2(\mu)}\leq C\,\|\RR_\mu^{d-1}\|_{L^2(\mu)\to L^2(\mu)},$$
where $C$ depends just on $d$ and $C(j)$.
\end{corollary}

Another related application of Theorem \ref{teomain22} appears in \cite{MoleroMourgoglouPuliattiTolsa}, by Molero, Mourgoglou, Puliatti, and the author,
where it is shown that the $L^2(\mu)$ boundedness of the singular integral operator $T_\mu$ whose kernel is given the gradient of the fundamental solution of
an elliptic operator in divergence form with double mean oscillation coefficients is equivalent to the $L^2(\mu)$ boundedness of the $(d-1)$-Riesz transform.

A possible approach to try to solve the David-Semmes problem for $n$-Riesz transforms with $1\leq n\leq d-1$ consists in studying the so-called reflectionless measures. For the purposes of this
paper, one can define a reflectionless measure for the Riesz transform $\RR^n$ as a measure $\mu$ in $\R^d$ having $n$-polynomial growth, such that $\RR^n_\mu$ is bounded in $L^2(\mu)$, and such that $\RR^n\mu$ is zero in a $\BMO(\mu)$ sense, i.e., $\RR^n\mu$ is constant $\mu$-a.e. 
Jaye and Nazarov have shown in \cite{JayeNazarov2} that if all the $n$-Ahlfors reflectionless measures for the Riesz transform are $n$-flat, i.e., of the form
$c\,\HH^n|_L$ for some $n$-plane $L$, then the David-Semmes problem for $\RR^n$ has a positive answer. However, up to now the structure of 
reflectionless measure for the Riesz transform $\RR^n$ is not well understood. It is not known if they are all $n$-flat even in the $n$-Ahlfors regular case with
$n=d-1$. 

\begin{question}
Let $n$ be an integer such that $2\leq n\leq d-1$. If $\mu$ is a reflectionless measure with polynomial $n$-growth (or even 
$n$-Ahlfors regular) for the $n$-Riesz transform, then is $\mu$
$n$-flat? 
\end{question}

Only in the case $n=1$ is the answer to this question known, and it is true, as shown by Jaye and Nazarov in \cite{JayeNazarov-Cauchy}. See also \cite{JayeNazarov1} 
for other rigidity properties of reflectionless measures, and \cite{MelnikovPoltoratskiVolberg} for a related result where Melnikov,  Poltoratski,
and Volberg
study the measures $\mu$ in $\R^d$ such that the Cauchy transform of $\mu$ vanishes $\mu$-a.e.\ in the principal value sense.


\section{The Painlev\'e problem for Lipschitz harmonic functions.} \label{sec5}

Given a compact set $E\subset \R^{d}$, one says that $E$ is removable for Lipschitz harmonic functions 
if for any open set $\Omega\supset E$, any function $f:\Omega\to\R$ which is Lipschitz in $\Omega$ and harmonic in
$\Omega\setminus E$
can be extended as a harmonic function to the whole $\Omega$. 
In the plane, $E$ is removable for bounded holomorphic functions  
if for any open set $\Omega\supset E$, any function $f:\Omega\to\R$ which is holomorphic and bounded in 
$\Omega\setminus E$
can be extended as a holomorphic function to the whole $\Omega$.
The Painlev\'e problem (for Lipschitz harmonic functions or holomorphic functions) consists in obtaining a metric-geometric description of removable sets for this class of functions.
To study this question and other related problems on approximation by Lipschitz harmonic and holomorphic functions, it is useful to introduce the Lipschitz harmonic capacity $\kappa$ and the analytic capacity $\gamma$. The Lipschitz harmonic capacity of $E$ is defined by
$$\kappa(E) = \sup |\langle \Delta f,1\rangle|,$$
where the supremum is taken over all Lipschitz functions $f:\R^{d}\to\R$ which are harmonic in $\R^{d}\setminus E$ and satisfy $\|\nabla f\|_\infty\leq1$, with $\Delta f$ understood in the sense of distributions. This capacity can be considered as the analogue of analytic capacity in higher dimensions. Recall that the analytic capacity of $E\subset\C$ equals
$$\gamma(E) = \sup |\langle \bar\partial f,1\rangle| = \sup |f'(\infty)|,$$
where the supremum is taken over all holomorphic functions in $\C\setminus E$ such that $\|f\|_{\infty,\C\setminus E}\leq 1$, assuming that $f=0$ on $E$ and that $\bar\partial f$ is taken in the sense of distributions, and  $f'(\infty) = \lim_{z\to\infty}z(f(z)-f(\infty))$.

Analytic capacity was introduced in 1947 by Ahlfors \cite{Ahlfors-bounded}, who showed that a compact set $E\subset\C$ is removable for bounded
analytic functions if and only if it has vanishing analytic capacity. Later, in the 1960's, Vitushkin \cite{Vitushkin} studied the uniform approximation in compact sets of holomorphic functions by rational functions and polynomials and he obtained criteria for approximation in terms of some estimates involving analytic capacity. On the other hand, Lipschitz harmonic capacity was introduced by Paramonov \cite{Paramonov} in 1990 in a work about approximation for
Lipschitz harmonic functions. See also \cite{MattilaParamonov}, where it is shown that $E$ is removable for Lipschitz harmonic functions if and only if $\kappa(E)=0$.

  The author \cite{TolsaSemiadditivity} (in the case $d=2$) and Volberg \cite{Volberg} (when $d>2$) showed  the following.
  
\begin{theorem}\label{teosem}  
Let $E\subset\R^{d}$ be compact. Then
$$\kappa(E) \approx \sup_{\mu\in A(E)} \mu(E),$$
where $A(E)$ is the family of all Radon measures $\mu$ supported on $E$ satisfying the $(d-1)$-polynomial growth condition \rf{eqgrow01} with constant $C_0=1$
 such that $\|\RR_\mu^{d-1}\|_{L^2(\mu)\to L^2(\mu)}\leq 1$. 
\end{theorem}

The result in the case $d=2$ also holds for analytic capacity. From this theorem it follows that the capacity $\kappa$ is countably
semiadditive. That is, there exists some $C>1$ such that
$$\kappa\Big(\bigcup_i E_i\Big) \leq C\,\sum_i\kappa(E_i).$$
The same result holds for analytic capacity.

For a compact set $E\subset \R^{d}$ with $\HH^{d-1}(E)<\infty$,  
any measure $\mu$ supported on $E$ with the polynomial growth condition \rf{eqgrow01} is of the
form $\mu= g\,\HH^{d-1}|_E$, for some $g\in L^\infty(\HH^{d-1}|_E)$. Together with Theorems  \ref{teosem} and \ref{teoriesz2}, this yields \cite{NazarovTolsaVolberg2}:

\begin{theorem}\label{teomain3}
Let $E\subset\R^{d}$ be compact with $\HH^{d-1}(E)<\infty$. 
Then $E$ is non-removable  
for Lipschitz harmonic functions if and only if there exists a $(d-1)$-rectifiable subset $F\subset E$ with $\HH^{d-1}(F)>0$.
\end{theorem}

Equivalently, given a compact set $E\subset\R^{d}$ with $\HH^{d-1}(E)<\infty$, 
$E$ is removable for Lipschitz harmonic functions if and only if $E$ is purely $(d-1)$-unrectifiable, i.e., $E$ does not contain any $(d-1)$-rectifiable subset with positive $\HH^{d-1}$ measure.

We remark that the case $d=2$ of Theorem \ref{teomain3} is due to David and Mattila \cite{DavidMattila}. Later David proved the analogous result for
bounded holomorphic functions 
in his celebrated work on Vitushkin's conjecture
\cite{David-vitus}. In both works  \cite{DavidMattila} and \cite{David-vitus}, the arguments depend strongly on the connection between the Cauchy kernel and Menger curvature.

The characterization of non-removability in terms of rectifiability obtained in Theorem \ref{teomain3} does not hold for sets with non-$\sigma$-finite measure $\HH^{d-1}$, although Theorem \ref{teosem} is still valid for these sets. 
In this case, combining Theorems \ref{teosem} and \ref{teomain22}, we deduce that a compact set in $\R^d$ is non-removable for Lipschitz harmonic functions if and only if it supports a non-trivial measure $\mu$ satisfying the growth condition \rf{eqgrow01} and such that
$$\int_0\int_0^\infty \beta_{2,\mu}^{d-1}(x,r)^2\,\Theta_\mu^{d-1}(x,r)\,\frac{dr}r\,d\mu(x)<\infty.$$
In fact, we deduce the following more precise result:

\begin{corollary}\label{coro2}
Let $E\subset\R^{d}$ be compact. Then
$$\kappa(E)\approx \sup\mu(E),$$
where the supremum is taken over all Radon measures $\mu$ supported on $E$ such that
$$
\mu(B(x,r))\leq r^{d-1}\quad \mbox{ for all $x\in\supp\mu$ and all $r>0$}
$$
and
$$
\int\!\!\int_0^\infty \beta_{2,\mu}^{d-1}(x,r)^2\,\Theta_\mu^{d-1}(x,r)\,\frac{dr}r\,d\mu(x)\leq \mu(E).
$$
\end{corollary}

To derive this corollary, remark that if $\mu$ satisfies the conditions above, by Chebyshev one deduces that there is a big piece
$F\subset E\cap\supp\mu$, with $\mu(F)\approx\mu(E)$, such that the measure $\wt\mu = \mu|_F$ satisfies \rf{eqbetawolff2}, and so
$\RR_{\wt \mu}^{d-1}$ is bounded in $L^2(\wt \mu)$. Hence $\kappa(E)\gtrsim \mu(F)\approx\mu(E)$. The converse
direction of the corollary is a straightforward consequence of Theorems \ref{teosem} and  \ref{teomain22}.

Since the conditions above on the measure $\mu$  are stable by bilipschitz maps, it follows that if $\vphi:\R^{d}
\to\R^{d}$ is bilipschitz, then
$$\kappa(E)\approx \kappa(\vphi(E))\quad\mbox{ for any compact set $E\subset\R^{d}$,}$$
with the implicit constant just depending on the bilipschitz constant of $\vphi$ and the ambient dimension. The same result had been obtained previously in \cite{Tolsa-bilip} in the case $d=2$, relying again on the connection between the Cauchy kernel and Menger curvature.

Another suggestive characterization of the capacity $\kappa(E)$ can be given in terms of the following non-linear potential, 
which we call the Jones-Wolff potential of $\mu$:
$$U_\mu(x) = \sup_{r>0} \Theta_\mu^{d-1}(x,r) + \left(\int_0^\infty \beta_{2,\mu}^{d-1}(x,r)^2\,\Theta_\mu^{d-1}(x,r)\,\frac{dr}r\right)^{1/2}.$$

\begin{corollary}\label{coro3}
Let $E\subset\R^{d}$ be compact. Then
$$\kappa(E) \approx \sup\{\mu(E):\,U_\mu(x)\leq 1\,\forall x\in E\}.$$
\end{corollary}

An immediate consequence of this result is that $E$ is non-removable for Lipschitz harmonic functions if and only if
it supports a non-zero measure $\mu$ such that $U_\mu(x)\leq1$ for all $x\in E$. 

The characterization of the capacity $\kappa$ and of removable sets for Lipschitz harmonic functions 
in terms of a metric-geometric potential such as $U_\mu$ should be considered as an analogue of the characterization
of analytic capacity and of removable sets for bounded analytic functions
in terms of curvature of measures \cite{TolsaSemiadditivity}. So one can think of the results stated in Corollaries
\ref{coro2} and \ref{coro3} as possible solutions of the Painlevé problem for Lipschitz harmonic functions.

Despite the advances described above, there are still interesting open questions related to the Painlev\'e problem for bounded holomorphic and Lipschitz harmonic functions. One of them deals with the connection between Favard length and analytic capacity. 
In the 1960's, Vitushkin conjectured that a compact set in the plane is non-removable for bounded holomorphic functions  if and only if its orthogonal projections have positive length in a set of directions of positive measure, or in other words, if and only if it has positive Favard length. 
The Favard length of a Borel set $E\subset\C$ equals
\begin{equation} \label{favard}
{\rm Fav}(E) =\int_0^\pi \HH^1(P_\theta(E))\,d\theta,
\end{equation}
where, for $\theta\in [0,\pi)$,  $P_\theta$ denotes the orthogonal
projection onto the line $L_\theta := \{re^{i\theta}:r\in\R\}$. 
For sets with finite length, Vitushkin's conjecture is true, because of David's theorem \cite{David-vitus}, since having zero Favard length for these sets is equivalent to being purely$1$-unrectifiable.
However, in 1986, Mattila \cite{Mattila-noVitushkin} showed that Vitushkin's conjecture is false for sets with non-$\sigma$-finite length. 
Indeed, he proved that the 
property of having zero Favard length is not invariant under conformal mappings while removability for bounded
analytic functions remains invariant. Mattila's result didn't tell which implication in the above conjecture
is false. This question was partially solved in 1988, when Jones and Murai \cite{JonesMurai} 
constructed a set with zero Favard length and positive analytic 
capacity. However, the following is still open:

\begin{question}
If $E$ is an arbitrary compact set in the plane with positive Favard length, then is its analytic capacity positive? 
\end{question}

Until recently, an important obstacle to solve the question above was the lack of a quantitative version of the Besicovitch projection theorem. A first step in this direction was obtained by Orponen \cite{Orponen-plentyprojections}, who proved the following result, solving a question put by David and Semmes in \cite{DavidSemmes2}.

\begin{theorem}
Let $E\subset \R^d$ be a compact $n$-Ahlfors regular set which has plenty of big projections (PBP). Then this is uniformly $n$-rectifiable (more precisely, it has big pieces of Lipschitz graphs).
\end{theorem}

 For simplicity, we only define the PBP property in the case $n=1$ in the plane: one says that $E$ has plenty of big projections if there exists $\delta>0$ such that for all $x\in E$ and $0<r\leq \diam(E)$ there exists
some $\theta_0\in [0,\pi)$ and some $\delta>0$ such that 
$$\HH^1(P_\theta(E\cap B(x,r)))\geq \delta\,r\quad \mbox{ for all $\theta\in [\theta_0-\delta,\theta_0+\delta].$}$$
Later on, relying on Orponen's result, D\k adrowski and Villa \cite{DabrowskiVilla} showed the following. Again we only announce the $1$-dimensional case in the plane.

\begin{theorem}
Let $E\subset\C$ be a compact set having PBP. Then $$\gamma(E)\gtrsim \diam(E),$$
with the implicit constant depending on the constants in the PBP property.
\end{theorem}

For other related results, see also \cite{OrponenMartikainen}, \cite{ChangTolsa}, and \cite{Tasso-projections}.

The problem of obtaining a quantitative version of the Besicovitch projection theorem was solved 
 by D\k abrowski in 2024 in \cite{Dabrowski-Favard}. In this work he proved the following remarkable result.

\begin{theorem}
Let $0<\kappa<1$ and $A\geq1$. Suppose that $E\subset\R^2$ is $1$-Ahlfors regular with constant $A$ and that ${\rm Fav}(E) \geq \kappa\,\HH^1(E)$.
Then there exists a Lipschitz graph $\Gamma$ with ${\rm Lip}(\Gamma)\leq C(\kappa,A)$ such that
$\HH^1(E\cap\Gamma)\geq c(\kappa,A)\,\HH^1(E)$, with $c(\kappa,A)>0$.
\end{theorem}

However, for the moment it is not clear if one can apply the techniques and ideas from D\k abrowski's work to make further progress in the solution of the open direction of Vitushkin's conjecture.


\section{Harmonic measure.} \label{sec6}

Let $\Omega\subset\R^{d}$ be an open set, and suppose for simplicity that it is bounded and Wiener regular, so that the Dirichlet problem for the Laplacian
is solvable for continuous functions on $\partial\Omega$. That is, for any function $f\in C(\partial\Omega)$ there exists a function
$u_f\in C(\overline \Omega)$ which is harmonic in $\Omega$ and such that $u_f|_{\partial\Omega} = f$. 
By the maximum principle, for any $x_0\in\Omega$, the functional $C(\pom)\ni f\mapsto u_f(x_0)\in \R$ is bounded with norm $1$ and it is positive.
Then, by the maximum principle and the Riesz representation theorem there exists a Radon probability measure $\omega^{x_0}$ supported on $\partial\Omega$ such that
\begin{equation}\label{eqharmonicmeasure}
u_f(x_0) = \int_{\partial\Omega} f\,d\omega^{x_0}
\end{equation}
for any $f\in C(\partial\Omega)$.
The measure $\omega^{x_0}$ is the harmonic measure for $\Omega$ with pole at $x_0$. We may drop the superscript $x_0$ in the notation if the choice of the pole is not relevant.
The notion of harmonic measure can also be introduced for domains which are non-Wiener regular via the Perron method and also for unbounded domains, but we skip these cases. For more information about this topic, we refer to the monograph \cite{PratsTolsa-llibre}.

Understanding the metric and geometric properties of harmonic measure is important for the study of the Dirichlet problem for the Laplace equation and also in complex analysis in the planar case $d=2$. In fact, in the plane, the behavior of conformal mappings at the boundary of $\Omega$ is strongly connected to the properties of the harmonic measure for $\Omega$. An old and natural question in classical analysis and PDE's  is to understand the relationship 
between harmonic and the Hausdorff measure $\HH^{d-1}$ (or other Hausdorff measures) on $\partial\Omega$. For example, an old theorem
of the Riesz brothers \cite{RieszRiesz} asserts that if $\Omega$ is a simply connected domain in the plane with $\HH^1(\partial\Omega)<\infty$, then the harmonic measure for $\Omega$ is mutually absolutely continuous with the arc-length measure $\HH^1$ on $\partial\Omega$.

A fundamental result in connection with harmonic measure is Dahlberg's theorem \cite{Dahlberg-Harmonicmeasure}. Roughly speaking, this asserts that for these domains 
harmonic measure and surface measure are mutually absolutely continuous in a quantitative way. The precise statement is the following.

\begin{theorem}\label{teodahlberg}
Let $\Omega\subset\R^{d}$ be a bounded Lipschitz domain and denote by $\sigma$ the surface measure
in $\pom$. Let $B$ be a ball centered in $\pom$ and let $x_0\in\Omega$ be such that $\dist(x_0,\frac32B\cap\pom)\geq C_1^{-1}\rad(B)$.
Then the following holds:
\begin{itemize}
\item[(a)] The harmonic measure $\omega^{x_0}$ and $\sigma$ are mutually absolutely continuous.
\item[(b)] We have
\begin{equation}\label{eqdahl}
\left(\avint_{B\cap\pom} \left(\frac{d\omega^{x_0}}{d\sigma}\right)^2\,d\sigma\right)^{1/2} \leq C\,\avint_{B\cap\pom} \frac{d\omega^{x_0}}{d\sigma}\,d\sigma = C\,\frac{\omega^{x_0}(B)}{\sigma(B)}.
\end{equation}
 Further, $\omega^{x_0}$ is a doubling measure. The constant $C$ above and the doubling constant depend only on $d$, the Lipschitz character of $\Omega$, and $C_1$.
\end{itemize}
\end{theorem}

In the statement (b), when we assert that $\omega^{x_0}$ is a doubling measure, we mean that there exists some constant $C$ such that
$$\omega^{x_0}(B(x,2r))\leq C\,\omega^{x_0}(B(x,r))\quad \mbox{ for all $x\in\pom$ and $0<r<\diam(\pom)$.}$$
Observe that \rf{eqdahl} is an $L^2$-reverse H\"older inequality.

There are many results which show that rectifiability plays a key role in the connection between harmonic measure and surface measure. In most of these results, Dahlberg's theorem is a basic ingredient. This is the case, for example, in the generalization of McMillan's theorem
to higher dimensions. To state it, we need to introduce the notion of cone point.

Given an open set $\Omega\subset \R^{d}$, a point $x\in\pom$ is called a cone point (for $\Omega$) if there exist a unit vector $v\in\R^{d}$, 
 $\alpha\in(0,1)$, and $r>0$ such that
$$X^+(x,v,\alpha,r):= \big\{y\in B(x,r):(y-x)\cdot v> \alpha\,|y-x|\big\}\subset \Omega.$$
Then we have:

\begin{theorem}\label{teoconepoints}
Let $\Omega\subset \R^{d}$ be a bounded domain, let $x_0\in\Omega$, and let $K\subset\pom$ be the subset of cone points for $\Omega$.
Then $\HH^{d-1}|_K \ll \omega_\Omega^{x_0}|_K$, that is, $\HH^{d-1}$ is absolutely continuous with respect to $\omega_\Omega^{x_0}$ on $K$.

Suppose moreover that the following holds: there exists some $c>0$ such that
\begin{equation}\label{eqbigcomp43}
\HH^{d-1}_\infty(B(\xi,r))\setminus \Omega) \geq c\,r^{d-1}\quad \mbox{ for all $\xi\in \pom$ and $r>0$.}
\end{equation}
Then $\omega_\Omega^{x_0}|_K\ll \HH^{d-1}|_K$, that is, $\omega_\Omega^{x_0}$  is absolutely continuous with respect to $\HH^{d-1}$ on $K$.
\end{theorem}

For simply connected domains in the plane, this theorem is a well-known result of McMillan \cite{McMillan}.
The extension to higher dimensions announced above
is due to Akman, Azzam, and Mourgoglou \cite{AkmanAzzamMourgoglou}. 
We recall that, for any open set $\Omega\subset \R^{d}$, the subset $K\subset\pom$  of cone points for $\Omega$ is $(d-1)$-rectifiable. This 
is a standard fact whose proof can be found in \cite[Chapter 15]{Mattila}, for example.

We will see in Section \ref{sec7} other results of more quantitative nature about the relationship between harmonic measure and surface measure.
Before turning to this type of results, we will deal with two free boundary problems for harmonic measure which have been solved by using the rectifiability criteria in terms of the Riesz transform described in Section \ref{sec4}.
The first one is a one-phase free boundary result which solves a question raised by Bishop in 1992 \cite{BishopQuestions}:

\begin{theorem}\label{teo-rectif}
Let $\Omega\subset\R^{d}$ be an open set and let $y_0\in\Omega$.
Suppose that there exists a set $E\subset\partial\Omega$ such that $0<\HH^{d-1}(E)<\infty$ 
and that the harmonic measure $\omega_\Omega^{y_0}|_E$ is mutually absolutely continuous with $\HH^{d-1}|_{E}$. Then
 $E$ is $(d-1)$-rectifiable.
\end{theorem}

This theorem was proved by Azzam, Hofmann, Martell, Mayboroda, Mourgoglou, the author, and Volberg \cite{AHMMMTV}. 
Notice that both the Riesz brothers theorem and McMillan's Theorem \ref{teoconepoints} assert that rectifiability plus other conditions imply the mutual absolute continuity between harmonic measure and the Hausdorff measure $\HH^{d-1}$. The theorem above deals with the converse direction:
the mutual absolute continuity of those measures implies rectifiability, for any arbitrary open set $\Omega$.

The connection between the Riesz transform and harmonic measure stems from the identity
$$G(x,y_0) = \EE(x-y_0) - \int_\pom \EE(x-z)\,d\omega^{y_0}(z),\quad \mbox{ for $x,y_0\in \Omega$,}$$
where $G$ is the Green function for $\Omega$. Differentiating on $x$, we get
\begin{equation}\label{eqgreenRiesz}
\nabla_x G(x,y_0) = \nabla \EE(x-y_0) - c_d\,\RR^{d-1}\omega^{y_0}(x),
\end{equation}
where $c_d$ is some absolute constant. In case $\Omega$ satisfies the CDC, from this identity and some standard estimates relating the Green function and harmonic measure, one deduces that
$$\RR^{d-1}_*\omega^{y_0}(x) \lesssim \sup_{r>0}\frac{\omega^{y_0}(B(x,r))}{r^{d-1}}.$$
In case the CDC does not hold, this estimate may fail but quite often one can get some weaker form of this inequality which may suffice for applications.

To prove Theorem \ref{teo-rectif}, in \cite{AHMMMTV} one uses the identity \rf{eqgreenRiesz}, the relationship between the Green function and harmonic measure, and the mutual absolute continuity between $\omega^{y_0}|_E$ and $\HH^{d-1}|_E$, to prove that 
$\RR^{d-1}_*\omega^{y_0}(x)<\infty$ for $\omega^{y_0}$-a.e.\  $x\in E$. Then, by Corollary \ref{coromaximal}, one deduces that 
$E$ is $(d-1)$-rectifiable.

Next we turn to the aforementioned two-phase free boundary problem for harmonic measure. The following theorem by Azzam, Mourgoglou, the author, and Volberg \cite{AzzamMourgoglouTolsaVolberg} solves a question posed again by Bishop in \cite{BishopQuestions}. The planar case $d=2$ had 
been proved previously by Bishop himself \cite{Bishop-Arkiv}.

\begin{theorem}\label{teo2phase}
For $d\geq 3$, let $\Omega_1,\Omega_2\subset \R^{d}$ be two disjoint domains and denote by 
$\omega_1$ and $\omega_2$ their respective harmonic measures. Let $E\subset \partial \Omega_1\cap
\partial \Omega_2$ 
 be a Borel set such that $\omega_{1}|_E$ and $\omega_{2}|_E$ are mutually absolutely continuous. Then $E$ contains a $(d-1)$-rectifiable subset $F$ with $\omega_1(E\setminus F)=0$ such that
$\omega_1|_F$ and $\omega_2|_F$ are mutually absolutely continuous with $\HH^{d-1}|_F$. 
\end{theorem}

We remark that under the additional assumption that $\Omega$ satisfies the capacity density condition, the theorem above had been proved previously by Azzam, Mourgoglou, and the author in \cite{AzzamMourgoglouTolsa-2phase}. Under this additional assumption, one obtains the existence
of tangents for $\Omega_1$ and $\Omega_2$ $\omega_i$-a.e.\ in $E$.  For other results of more quantitative nature about the two phase problem for harmonic measure, see for example \cite{KenigToro-2phase}, \cite{BortzHofmann},
\cite{AzzamMourgoglouTolsa-2phasequant}, \cite{PratsTolsa2}, \cite{TolsaToro-2phase}.

By a previous result of Kenig, Preiss, and Toro \cite{KenigPreissToro}, when $\Omega_1$ and $\Omega_2$ are NTA domains (see Section \ref{sec6} for the precise definition), it was already known that the mutual absolute continuity of the harmonic measures $\omega_1$ and $\omega_2$ implies that these measures are concentrated in
a $(d-1)$-dimensional set. 
The main difficulty in proving Theorem \ref{teo2phase} consists in showing the existence of a subset $F\subset E$ with $\omega_1(E\setminus F)=0$ such that
$\omega_1|_F$ and $\omega_2|_F$ are mutually absolutely continuous with $\HH^{d-1}|_F$, since then the $(d-1)$-rectifiability of $F$
follows from Theorem~\ref{teo-rectif}. By standard arguments, it suffices to show that
$$0<\theta^{*,d-1}(x,\omega_1):= \limsup_{r\to0}\frac{\omega_1(B(x,r))}{(2r)^{d-1}}<\infty\quad \mbox{ for $\omega_1$-a.e.\ $x\in E$.}$$
The main difficulty is in proving the left inequality. To this end, the
main tool is  Theorem \ref{teoGT}.
Indeed, suppose that there exists some 
subset $F\subset\pom_1\cap\pom_2$ with positive measure $\omega_1$ such that $\theta^{*,d-1}(x,\omega_1)=0$ for al $x\in F$. 
Then, an application of Theorem \ref{teoGT} in a suitable small ball $B$ centered at a density point of $F$ would lead to a contradiction, since this ensures that a non-negligible portion $G_0\subset F$ is supported on a uniformly rectifiable set $\Gamma$, which implies that $\theta^{*,d-1}(x,\omega_1)>0$
a.e.\ on $G$. To guarantee that the assumptions of the theorem are satisfied, one needs to apply the Alt-Caffarelli-Friedman monotonicity formula and a delicate blowup argument, in a similar fashion to \cite{KenigPreissToro}, together with the connection between the Riesz transform and harmonic measure that stems from \rf{eqgreenRiesz}.

To conclude this section, we wish to recall an important open question about the dimension of harmonic measure.
For any Borel measure $\nu$
in $\R^{d}$, its (Hausdorff) dimension, denoted by $\dim_H\nu$, is defined by
$$\dim_H\nu = \inf\,\{\dim_H F: F\subset \R^{d} \text{ Borel, } \nu(F^c)=0\},$$
where $\dim_H F$ stands for the Hausdorff dimension of $F$.
A celebrated theorem of Jones and Wolff \cite{JonesWolff} asserts that, for any arbitrary domain in the plane, $\dim_H\omega\leq1$.
In the higher dimensional case $d>2$, the situation is more complicated. 
On the one hand, Bourgain \cite{Bourgain-dimension} proved in 1987 that there exists some constant $\ve_d>0$ just depending on $d$ such that
$\dim_H\omega\leq d-\ve_d$ for any $\Omega\subset\R^{d}$.
A natural guess would be that one could take $\ve_d=1$, so that
$\dim_H\omega\leq d-1$ for any domain of $\R^{d}$, analogously to what happens in the plane. However, this was disproved by Wolff in 
\cite{Wolff-counterexample}, where
he managed to construct a snowflake type domain $\Omega\subset\R^{d}$ satisfying $\dim_H\omega>d-1$, for $d\geq3$.
A difficult open question in the area is the following:

\begin{question}
Which is the optimal value of the constant $\ve_d$? A less demanding but related question is to determine its asymptotic behavior as the dimension $d$ grows. 
\end{question}

Bishop \cite{BishopQuestions} and Jones \cite{Jones-scaling} conjecture that $\ve_d=\frac1{d-1}$. However, at present there is not much evidence supporting this conjecture.


\section{Boundary value problems for harmonic functions in rough domains.}\label{sec7}

Let $\Omega\subset\R^d$ be a bounded open set such that $\HH^{d-1}(\pom)<\infty$. Denote $\sigma=\HH^{d-1}|_\pom$ and $L^p(\pom)= L^p(\sigma)$ for $1\leq p \leq\infty$.
Recall that the Dirichlet problem (for the Laplace equation) with boundary data $f:\partial \Omega\to\R$, $f\in C(\partial\Omega)$ consists in finding the function $u\in
C^2(\Omega)\cap C(\overline \Omega)$ such that 
$$\left\{\begin{array}{ll}\Delta u= 0 & \quad\text{ in $\Omega$,}\\
u= f  & \quad\text{ in $\partial\Omega$.}
\end{array}
\right.
$$
On the other hand, the Neumann problem with data $g:\partial \Omega\to\R$, with $g\in L^2(\partial\Omega)$ such that $\int_{\partial\Omega} g\,d\sigma=0$ (we write $g\in L^2_0(\partial\Omega)$),
consists in finding the function $u\in C^2(\Omega)$ such that
$$\left\{\begin{array}{ll}\Delta u= 0 & \quad\text{ in $\Omega$,}\\
\partial_\nu u= g  & \quad\text{ in $\partial\Omega$,}
\end{array}
\right.
$$
where $\partial_\nu u$ stands for the outer normal derivative of $u$ in $\pom$. This normal derivative can be defined as a non-tangential limit if it exists, or in a weak sense, by the identity
\begin{equation}\label{eqweak93}
\int_\Omega \nabla u\,\nabla \vphi\,dx = \int_\pom g\,\vphi\,d\sigma \quad \mbox{ for all $\vphi\in C_c^\infty(\R^d)$.}
\end{equation}

For $\Omega$ as above, a parameter $a>0$, a function $u:\Omega\to\R$, and $\xi\in \pom$, we consider the non-tangential maximal operator
$$\cN(u)(\xi) = \sup_{y\in \Gamma_a(\xi)} |u(y)|,$$
where 
$\Gamma_a(\xi) = \{y\in\Omega: |\xi-y|< (1+a)\,\dist(y,\partial\Omega)\}.$
Notice that $\Gamma_a(\xi)$ is a non-tangential region with vertex at $\xi$.
We say that
the Dirichlet problem (for the Laplace equation) is $L^p$-solvable if the solution $u$ of the Dirichlet problem with boundary data $f$  satisfies
$$\|\cN(u)\|_{L^p(\pom)} \leq C\,\|f\|_{L^p(\pom)}
\quad \mbox{ for all Dirichlet data $f\in C(\partial \Omega)$.}$$
The regularity problem is $L^p$-solvable if
$$\|\cN(\nabla u)\|_{L^p(\partial\Omega)} \leq C\,\|\nabla_t f\|_{L^p(\partial\Omega)}
\quad \mbox{ for all Dirichlet data $f\in {\rm Lip}(\partial \Omega)$,}$$
where $\nabla_t f$ is the tangential derivative of $f$ in $\pom$. This tangential derivative can be defined pointwise $\sigma$-a.e.\ when $\pom$
is rectifiable, or alternatively, one can use the framework of the Hajlasz-Sobolev spaces \cite{Hajlasz} to define this.
On the other hand, the Neumann problem is $L^p$-solvable if
$$\|\cN(\nabla u)\|_{L^p(\partial\Omega)} \leq C\,\|\partial_\nu u\|_{L^p(\partial\Omega)}
\quad \mbox{ for all Neumann data $g=\partial_\nu u\in L_0^2(\partial \Omega)$.}$$
For short, we say that $(D_p)$, $(R_p)$, $(N_p)$ are solvable, respectively. It is well known that if $\pom$ is Ahlfors regular, then the solvability of these problems does not depend on the aperture $a$ of the non-tangential region used to define $\cN$ (although the constants $C$ above may depend on it).

To study the solvability of the Dirichlet, regularity, and Neumann problems one needs to assume quantitative geometric conditions on the 
underlying domain. Recall that one says that a set $\Omega$ satisfies the corkscrew condition, or that $\Omega$ is a corkscrew set, if
there exists some constant $c_1>0$ such that for all
$x\in\pom$ and all $r\in(0, 2\diam(\Omega))$ there exists a ball $B\subset B(x,r)\cap\Omega$ such that
$\rad(B)\geq c_1\,r$. On the other hand, $\Omega$ is called two-sided corkscrew if both $\Omega$ and $\R^d\setminus \overline{\Omega}$ are corkscrew sets.

Given two points $x,x' \in \Omega$, and a pair of numbers $M,N\geq1$, 
an $(M,N)$-{Harnack chain connecting $x$ to $x'$},  is a chain of
open balls
$B_1,\dots,B_N \subset \Omega$, 
with $x\in B_1,\, x'\in B_N,$ $B_k\cap B_{k+1}\neq \varnothing$
and $M^{-1}\diam (B_k) \leq \dist (B_k,\partial\Omega)\leq M\diam (B_k)$ for all $k$.
For $C_2\geq1$, we say that $\Omega$ satisfies the {\it $C_2$-Harnack chain condition} if
 for any two points $x,x'\in\Omega$,
there is an $(M,N)$-Harnack chain connecting them, with $M\leq C_2$ and $N$ such that
$$N\leq C_2\,\bigg(1+\log^+\frac{|x-x'|}{\min(\delta_\Omega(x),\delta_\Omega(x'))}\bigg),$$
where $\delta_{\Omega}(x)=\dist(x, \pom)$.

 A domain $\Omega\subset\R^d$ is called uniform if it is a corkscrew domains and it satisfies the Harnack chain condition.
Following \cite{JerisonKenig}, $\Omega$ is called NTA ({\it non-tangentially accessible})  if it is uniform and two-sided corkscrew.  One says that set $\Omega \subset \R^{d}$
 is a chord-arc domain if it is NTA and  $\pom$ is $(d-1)$-Ahlfors regular. 
  Additionally, if a domain $\Omega$ and its exterior $\R^{d} \setminus \overline\Omega$  are chord-arc, then  $\Omega$ is called two-sided chord-arc. 
  
\subsection{The Dirichlet problem.}

The following theorem characterizes the $L^p$ solvability of the Dirichlet problem in terms of a local reverse $L^{p'}$ H\"older inequality for harmonic measure.

\begin{theorem}\label{teo9.2}
Let $\Omega\subset\R^d$ be a domain with bounded $(d-1)$-Ahlfors regular boundary and set $\sigma=\HH^{d-1}|_\pom$. Given $x \in\Omega$, denote by $\omega^x$ the
harmonic measure for $\Omega$ with pole at $x$.
For $p\in (1,\infty)$, the following are equivalent:
\begin{itemize}
\item[(a)] $(D_{p})$ is solvable for $\Omega$.

\item[(b)] The harmonic measure $\omega$ is absolutely continuous with respect to $\sigma$ and
for every ball $B$ centered in $\pom$ and for all $x\in \Omega\cap 3B\setminus 2B$ with $\diam(B)\leq 2\diam(\pom)$, it holds that
\begin{equation}\label{equf88noDuality}
\left(\avint_{B} \left(\frac{d\omega^{x}}{d\sigma}\right)^{p'}\,d\sigma\right)^{1/p'} \lesssim \sigma(B)^{-1},
\end{equation}
where $p'=p/(p-1)$.
\end{itemize}
\end{theorem}

This theorem is proven in \cite{MourgoglouTolsa-regu}, but quite likely it was already known before (this is sure at least for the class of chord-arc domains). See, for example, 
\cite{Kenig} and \cite{Hofmann-sinica} for somewhat less precise results.
Remark that combining Dahlberg's Theorem \ref{teodahlberg} and the preceding result, one deduces the (well known) fact that the Dirichlet problem is solvable in $L^p$ in Lipschitz domains for $p\geq2$.

For an open set $\Omega\subset\R^d$ with $(d-1)$-Ahlfors regular boundary, one says that the harmonic measure belongs to local weak-$A_\infty$
if the condition (b) in Theorem \ref{teo9.2} holds for some $p>1$. In the case when $\Omega$ is a chord-arc domain, using the change of pole formula
and the fact that harmonic measure is doubling (see \cite{JerisonKenig}), it follows easily that the local weak-$A_\infty$ condition is equivalent to the fact that, for a fixed pole $x_0\in\Omega$, the harmonic measure $\omega^{x_0}$ is an $A_\infty$ weight with respect to the surface measure $\HH^{d-1}|_\pom$.

It was shown independently by David and Jerison \cite{DavidJerison} and Semmes \cite{Semmes} that if $\Omega$ is two-sided corkscrew, then its boundary $\pom$ is uniformly $(d-1)$-rectifiable, and moreover $\Omega$ has big pieces of Lipschitz subdomains. The last statement means that
for any ball $B(x,r)$ centered in $\pom$ with radius at most $\diam(\pom)$ there exists a Lipschitz subdomain $\Omega_{x,r}\subset\Omega$ (with constants uniform on $x$ and $r$) such that, for some fixed $\gamma>0$,
\begin{equation}\label{eqbigpieces}
\HH^{d-1}(\pom\cap \pom_{x,r}\cap B(x,r))\geq \gamma\,\HH^{d-1}(\pom\cap B(x,r)).
\end{equation}
Using this fact, it was shown in \cite{DavidJerison} and \cite{Semmes} that if $\Omega$ is chord-arc domain, then the harmonic measure satisfies the local weak-$A_\infty$ property, and so it belongs to $A_\infty$ in this case. Later on, 
 Bennewitz and Lewis \cite{BennewitzLewis} proved that the two-sided corkscrew condition (together with the $(d-1)$-Ahlfors regularity of the boundary) suffices for the local weak-$A_\infty$ condition to hold.
 
The full geometric characterization of the domains with $(d-1)$-Ahlfors regular boundaries for which the local weak-$A_\infty$ condition holds was obtained by
Azzam, Hofmann, Martell, Mourgoglou, and the author in \cite{AzzamHofmannMartellMourgoglouTolsa}. The precise result is the following.

\begin{theorem}\label{teoAHMMT} 
Let $\Omega\subset \R^d$, $d\ge 3$, be an open set with $(d-1)$-Ahlfors regular boundary satisfying the  
corkscrew condition.  Then the following are equivalent:
\begin{itemize}
\item[(a)] $\Omega$  has big pieces of 
chord-arc subdomains.
 
\item[(b)] Harmonic measure for $\Omega$ is in local weak-$A_\infty$. 

\end{itemize}
\end{theorem}

The notion of big pieces of chord-arc subdomains is analogous to the one of big pieces of Lipschitz subdomains introduced in \rf{eqbigpieces}: one only has to replace the assumption that the sets $\Omega_{x,r}$ are Lipschitz subdomains by the one of being chord-arc subdomains (with uniform constants too). We remark that in \cite{AzzamHofmannMartellMourgoglouTolsa} it is also shown that the conditions (a) and (b) above are equivalent to the fact that $\pom$ is uniformly $(d-1)$-rectifiable and $\Omega$ satisfies a suitable connectivity condition, namely, the so-called weak local John condition. The fact that
the local weak-$A_\infty$ condition implies uniform rectifiability had been shown previously in \cite{HLMN} and \cite{MourgoglouTolsa-reine}.
The proof of the implication (b) $\Rightarrow$ (a) in the theorem uses again the Alt-Caffarelli-Friedman monotonicity property together with some ideas originating from Aikawa and Hirata \cite{AikawaHirata} and other techniques from quantitative rectifiability. 

Many other results connecting the geometry of domains and PDE type results have been obtained in recent years. For example, the next result 
characterizes the uniform rectifiability of $\pom$.

\begin{theorem} 
Let $\Omega\subset\R^{d}$, $d\geq2$, be a corkscrew domain with $(d-1)$-Ahlfors regular boundary.  Then the following are equivalent:
\begin{itemize}
\item[(a)] $\partial\Omega$ is uniformly $(d-1)$-rectifiable.

\item[(b)] There is $C>0$ such that for every bounded harmonic function $u$ on $\Omega$ and every ball $B$ centered at $\partial\Omega$,
\begin{equation*}
\int_B |\nabla u(x)|^2\,\dist(x,\partial\Omega)\,dx\leq C\,\|u\|^2_{L^\infty(\Omega)}\,\rad(B)^{d-1}.
\end{equation*}
\end{itemize}
\end{theorem}

The implication (a) $\Rightarrow$ (b) is due to Hofmann, Martell, and Mayboroda \cite{HofmannMartellMayboroda}, and the converse one to Garnett, Mourgoglou, and the author \cite{GarnettMourgoglouTolsa}.

\subsection{The regularity and Neumann problems.}

For Lipschitz domains, in 1984 Verchota \cite{Verchota} proved that the regularity problem is solvable in $L^2$, and later on, in 1987, Dahlberg and Kenig \cite{DahlbergKenig} proved that both the regularity and the Neumann problems are solvable in $L^p$ in the range $p\in (1,2+\ve)$. 
Previously, Fabes, Jodeit, and Rivi\`ere had shown in \cite{FabesJodeitRiviere} that in $C^1$ domains both problems are solvable for all $p\in(1,\infty)$. The arguments in all these works rely on the proof of the invertibility of the double layer potential in $L^p$.

An important open problem in the area of boundary value problems consists in finding a larger class of domains where the regularity and Neumann problems for the Laplace equation are solvable in $L^p$. In particular, in 1991 Kenig  asked the following question  \cite[Problem 3.2.2]{Kenig}: in a bounded chord-arc domain $\Omega \subset \R^{d}$ does there exist $p>1$ such that the regularity and Neumann problems for the Laplacian are solvable in $L^p$?

To deal with the regularity problem in chord-arc domains it is convenient to introduce a variant of the problem in terms of the Haj\l asz-Sobolev spaces.
Fix a metric space $\Sigma$ and a doubling measure $\sigma$ with $\supp\sigma=\Sigma$ (below we will take $\Sigma = \pom$ and $\sigma=\HH^{d-1}|_\pom$). For 
a Borel function $f:\Sigma\to\R$, we say that a non-negative Borel function $g:\Sigma \to \R$ is a Haj\l asz upper gradient of  $f$ if  
\begin{equation}\label{eq:H-Sobolev}
 |f(x)-f(y)| \leq |x-y| \,(g(x)+g(y))\quad \mbox{ for $\sigma$-a.e. $x, y \in \Sigma$.} 
\end{equation}
We denote the collection of all the Haj\l asz upper gradients of $f$ by $D(f)$.
Then for $p\geq1$, we let $\dot{W}^{1,p}(\Sigma)$ be the space of Borel functions $f$ which have 
a Haj\l asz upper gradient in $L^p(\sigma)$ 
and we  define the semi-norm 
\begin{equation}\label{eqseminorm}
 \| f \|_{ \dot W^{1.p}(\Sigma)} = \inf_{g \in D(f)} \| g\|_{L^p(\Sigma)}
 \end{equation}

 We say that the problem $(\wt R)_p$ is solvable in $\Omega$ if the solution of the Dirichlet problem with boundary data $f\in {\rm Lip}(\partial \Omega)$ satisfies
$$\|\cN(\nabla u)\|_{L^p(\partial\Omega)} \leq C\,\| f\|_{\dot W^{1,p}(\pom)}.$$
The aforementioned question about the regularity problem was solved by Mourgoglou and the author in \cite{MourgoglouTolsa-regu}:

\begin{theorem}\label{teoregu}
Let $\Omega\subset\R^{d}$ be a bounded corkscrew domain with $(d-1)$-Ahlfors regular boundary.
For $p \in (1, \infty)$, the Dirichlet problem $(D_{p'})$ is solvable if and only if $(\wt R_p)$ is solvable. 
\end{theorem}

Since, for chord-arc domains, the Dirichlet problem is solvable in $(D_{p'})$ for $p'$ large enough, it follows that $(\wt R_p)$ is solvable in this type of domain for $p$ small enough. 
On the other hand, as shown in \cite{MourgoglouTolsa-regu}, in the absence of connectivity condition on $\pom$, the implication $(D_{p'})\,\Rightarrow\,(R_p)$ may fail. However, if $\pom$ supports a weak $p$-Poincar\'e inequality (see Section \ref{secprelim}), then $\|\nabla_t f\|_{L^p(\pom)}\approx \| f\|_{\dot W^{1,p}(\pom)}$ and thus 
the solvability of $(R_p)$ is equivalent to the solvability of $(\wt R_p)$. It has been shown in \cite{TapiolaTolsa} that if $\Omega$ is a two-sided chord-arc domain, then $\pom$ supports a weak $q$-Poincar\'e inequality for every $q\geq1$. So for two-sided chord-arc domains, $(R_p)$ is solvable for $p>1$ small enough.
We also remark that Azzam \cite{Azzam-poincare} has proved  that if $E$ is an $n$-Ahlfors regular set in $\R^d$ that supports a weak $n$-Poincar\'e inequality, then $E$ is uniformly $n$-rectifiable.

The  strategy for the proof of Theorem \ref{teoregu}  consists in constructing an ``almost harmonic extension'' $v$ of the boundary function $f \in \Lip(\pom)$  to $\Omega$, such that its distributional Laplacian $\Delta v$
satisfies an $L^p$-Carleson condition in $\Omega$ and such that, moreover,  its normal derivative $\partial_\nu v$ in $\pom$
is controlled by the Haj\l asz tangential gradient of $f$. 
The construction of this almost orthogonal extension is performed by means of a corona type decomposition of $\Omega$ in terms of mutually disjoint interior Lipschitz subdomains together with a ``buffer'' region. With the almost harmonic extension in hand, one is able to apply a duality argument to show that
Dirichlet solvability in $L^{p'}$ implies the existence the following  one-sided Rellich type inequality:
\begin{equation}\label{eqrellich}
\|\partial_\nu u\|_{L^p(\pom)} \leq C\,\| f\|_{\dot W^{1,p}(\pom)}.
\end{equation}
Then, using layer potential techniques, one deduces that $(\wt R_p)$ is solvable.

Regarding the Neumann problem, there has been little progress since the proof of the solvability for $p\in (1,2+\ve)$ in Lipschitz domains by Dahlberg and Kenig. In \cite{HofmannMitreaTaylor} Hofmann, Mitrea, and Taylor  proved that the Neumann problem is solvable in SKT domains for all $p\in(1,\infty)$.
An SKT domain is a chord-arc domain such that the outer unit normal belongs to $VMO$, which in particular, implies that $\pom$ is Reifenberg flat  \cite{HofmannMitreaTaylor}. That is, $\lim_{r\to0}\beta_{\infty,\pom}^{d-1}(x,r)=0$  for all $x\in\pom$, uniformly on $x$.

\begin{question}
Is there some $p>1$ small enough such that the Neumann problem is solvable in $L^p$ for chord-arc domains, with $p$ possibly depending on the chord-arc character?
\end{question}

The most general class of domains where the Neumann problem is known to be solvable in $L^p$ for some $p>1$ is the one of chord-arc domains such that the regularity problem is solvable in $L^p$ and have {\it very} big pieces
of Lipschitz subdomains, as shown by recently by Mourgoglou and the author \cite{MourgoglouTolsa-Neumann}. Recall that any chord-arc domain
$\Omega$ has big pieces of Lipschitz graph. That is, 
 for any ball $B(x,r)$ centered in $\pom$ with radius at most $\diam(\pom)$ there exists a Lipschitz subdomain domain $\Omega_{x,r}\subset\Omega$ (with constants uniform on $x$ and $r$) such that, for some fixed $\gamma\in (0,1)$,  
\rf{eqbigpieces} holds. One says that $\Omega$ has {\it very} big pieces of Lipschitz subdomains if \rf{eqbigpieces} holds for some $\gamma$ very close to $1$.
Of course, this definition depends on the particular choice of $\gamma$. 
\vv

To close this section, we want to remark that the area of boundary value problems for elliptic PDE's and its interaction with uniform rectifiability has been very active in the last years. Some directions of research that we have not mentioned above include the study of elliptic measures under the Dahlberg, Kenig, Pipher condition, with contributions such as the one of Hofmann, Martell, Mayboroda, Toro, and Zhao \cite{HMMTZ};
the one-phase and two-phase problems for elliptic measures associated to PDE's with VMO or double mean oscillation coefficients by Azzam and Mourgoglou \cite{AzzamMourgoglou-VMO}, and Merlo, Mourgoglou, and Puliatti \cite{MerloMourgoglouPuliatti};
the study of elliptic measures 
in domains with high codimension, with works by David, Engelstein, Feneuil, Mayboroda, and others, such as \cite{DavidFeneuilMayboroda-Dahlberg}, \cite{DavidFeneuilMayboroda2}, \cite{DavidEngelsteinMayboroda}; or the study of caloric measure
and parabolic uniform rectifiability in graph domains by Bortz, Hofmann, Martell, and Nystr\"om \cite{BortzHofmannMartellNystrom}. Other related works are \cite{MourgoglouPoggiTolsa} by Mourgoglou, Poggi, and the author and \cite{FeneuilLi} by Feneuil and Li, where some suitable Poisson problems are considered in connection with the regularity and Neumann problems, respectively. Also, in the work \cite{Gallegos-Rellich},  Gallegos studies the solvability of the regularity problem
in the end-point case $p=1$ and he shows just that the 
uniform rectifiability of $\pom$ suffices for the validity of the
one-sided Rellich inequality \rf{eqrellich} for $p=1$, and conversely, in the plane, the validity of such inequality for $p=1$ implies the uniform rectifiability of $\pom$.


\section*{Acknowledgments}

I am grateful to Mihalis Mourgoglou and Damian D\k{a}browski for some pertinent remarks on a preliminary version of this paper.

This work was supported by the European Research Council (ERC) under the European Union's Horizon 2020 research and innovation programme (grant agreement 101018680). Also partially supported by MICIU (Spain) under the grant PID2024-160507NB-I00 and the María de Maeztu Program for units of excellence (Spain) (CEX2020-001084-M).

\bibliographystyle{siamplain}

\bibliography{Llibres}
\end{document}